\documentclass[11pt,a4paper]{article}
\usepackage{amssymb}
\usepackage{amsmath}
\usepackage{latexsym}
\usepackage{amsthm}
\usepackage{todonotes}
\usepackage{comment}
\usepackage{xcolor}
\definecolor{myred}{rgb}{0.2,0,0}
\definecolor{myblue}{rgb}{0,0,0.6}
\definecolor{mygreen}{rgb}{0,0.2,0}

\newcommand{\ls}{\leqslant}
\newcommand{\gs}{\geqslant}

\newtheorem{thm}{Theorem}
\newtheorem{lem}[thm]{Lemma}
\newtheorem{cor}[thm]{Corollary}
\newtheorem{prop}[thm]{Proposition}
\theoremstyle{definition}

\usepackage{graphicx}
\usepackage{url}
\usepackage{xpatch}
\xpatchcmd{\proof}{\itshape}{\normalfont\proofnameformat}{}{}
\newcommand{\proofnameformat}{}

\usepackage{xcolor}
\usepackage{caption}
\usepackage{varwidth}

\newcommand{\captionbackgroundcolor}[1]{\colorlet{cpbgcol}{#1}}
\DeclareCaptionFont{black}{\color{black}}
\captionbackgroundcolor{black!3}

\newcommand\blfootnote[1]{%
  \begingroup
  \renewcommand\thefootnote{}\footnote{#1}%
  \addtocounter{footnote}{-1}%
  \endgroup
}
\makeatletter
\DeclareRobustCommand{\cev}[1]{%
  \mathpalette\do@cev{#1}%
}
\newcommand{\do@cev}[2]{%
  \fix@cev{#1}{+}%
  \reflectbox{$\m@th#1\vec{\reflectbox{$\fix@cev{#1}{-}\m@th#1#2\fix@cev{#1}{+}$}}$}%
  \fix@cev{#1}{-}%
}

\newcommand{\fix@cev}[2]{%
  \ifx#1\displaystyle
    \mkern#23mu
  \else
    \ifx#1\textstyle
      \mkern#23mu
    \else
      \ifx#1\scriptstyle
        \mkern#22mu
      \else
        \mkern#22mu
      \fi
    \fi
  \fi
}

\begin{document}

\renewcommand{\proofnameformat}{\bfseries}

\begin{center}
{\Large\textbf{On the behavior of the colored Jones polynomial of the figure-eight knot under modular transformations}}


\vspace{5mm}

\textbf{Christoph Aistleitner and Manuel Hauke}

\vspace{5mm}

{\footnotesize Graz University of Technology

Institute of Analysis and Number Theory

Steyrergasse 30, 8010 Graz, Austria

Email: \texttt{aistleitner@math.tugraz.at} and \texttt{hauke@math.tugraz.at}}
\end{center}

\begin{abstract}
The colored Jones polynomial $J_{K,N}$ is an important quantum knot invariant in low-dimensional topology. In his seminal paper on quantum modular forms, Zagier predicted the behavior of $J_{K,0}(e^{2 \pi i x})$ under the action of $\textup{SL}_2(\mathbb{Z})$ on $x \in \mathbb{Q}$. More precisely, Zagier made a prediction on the asymptotic value of the quotient $J_{K,0}(e^{2 \pi i \gamma(x)})/ J_{K,0}(e^{2 \pi i x})$ for fixed $\gamma  \in \textup{SL}_2(\mathbb{Z})$, as $x \to \infty$ along rationals with bounded denominator. In the case of the figure-eight knot $4_1$, which is the most accessible case, there is an explicit formula for $J_{4_1,0}(e^{2 \pi i x})$ as a sum of certain trigonometric products called Sudler products. By periodicity, the behavior of $J_{4_1,0}(e^{2 \pi i x})$ under the mapping $x \mapsto x+1$ is trivial. For the second generator of $\textup{SL}_2(\mathbb{Z})$, Zagier conjectured that with respect to the mapping $x \mapsto 1/x$, the quotient $h(x) = \log ( J_{4_1,0}(e^{2 \pi i x}) / J_{4_1,0}(e^{2 \pi i /x}))$ can be extended to a function on $\mathbb{R}$ that is continuous at all irrationals. This conjecture was recently established by Aistleitner and Borda in the case of all irrationals that have an unbounded sequence of partial quotients in their continued fraction expansion. In the present paper we prove Zagier's continuity conjecture in full generality. 
\end{abstract}

\section{Introduction and statement of results}

Among the topological invariants that are connected with a knot $K$ in $\mathbb{R}^3$, two of the most important ones are the colored Jones polynomial $J_{K,N}, ~N \geq 2,$ and the Kashaev invariant $\{\langle K \rangle\}_{N \geq 2}$. They are connected to each other via $\langle K \rangle = J_{K,N}(e^{2 \pi i / N})$, and are also related to the  Alexander polynomial, another important knot invariant. The volume conjecture, which is only solved in some special cases, relates the asymptotic behavior of the Kashaev invariant of a knot to the hyperbolic geometry of its complement, thereby suggesting that the Kashaev invariant and the colored Jones polynomial both encode information on the geometry of the knot complement. The volume conjecture has deep implications in mathematics and theoretical physics, including quantum gravity and topological quantum field theory. For more information on this general background and the volume conjecture, we refer to the monograph \cite{muyo} and the research papers \cite{betdr1, mura, muramura, witten, yokota}.\\

The Kashaev invariant can be extended to a function on the roots of unity by setting, for $\gcd(h,k) =1$, $\mathbf{J}_{K} (h/k) := J_{K,k}(e^{2 \pi i h/k})$. Zagier defined further $J_{K,0}(e^{2 \pi i h/k}) := J_{K,k}(e^{2 \pi i h/k})$ by backwards extrapolation (this is the function $J_{K,0}$ appearing in the abstract). The volume conjecture then predicts the size of the limit
$$
\lim_{N \to \infty} \frac{\log |\langle K \rangle_N|}{N} = \lim_{N \to \infty} \frac{\log |\mathbf{J}_{K}(1/N)|}{N}. 
$$
However, it seems that $\mathbf{J}_{K}$ has arithmetic properties that go far beyond this asymptotic relation. Concerning the behavior of $\mathbf{J}_{K}$ at other roots of unity, Zagier predicted an asymptotic formula for the quotient of $\mathbf{J}_{K}\left(\gamma x \right)$ and $\mathbf{J}_{K}(x)$ as $x \to \infty$  along rationals with bounded denominator, where $\gamma = \left(\begin{smallmatrix}a & b \\ c & d \end{smallmatrix}\right) \in \textup{SL}_2(\mathbb{Z})$ is fixed and acts on $\mathbb{Q}$ as $\gamma x = \frac{ax+b}{cx+d}$. This ``approximate modularity'' has been showcased in Zagier's paper \cite{ZA} on what he called ``quantum modular forms'', where the behavior of $\mathbf{J}_{K}$ under the action of $\textup{SL}_2(\mathbb{Z})$ is regarded as ``the most mysterious and in many ways the most interesting'' among the examples mentioned in the paper. Throughout the rest of this paper, we will only be concerned with the ``figure-eight knot'', written as $4_1$ in Alexander-Briggs notation. In many regards, this knot is the simplest non-trivial hyperbolic knot. For this particular knot, we have the explicit formula
\begin{equation} \label{q-series}
\mathbf{J}_{4_1} (x) = \sum_{n=0}^\infty |(1-e^{2 \pi i x}) (1-e^{2 \pi i 2 x}) \cdots (1-e^{2 \pi i n x})|^2, 
\end{equation}
for $x \in \mathbb{Q}$ (note that this actually is a finite sum, since all but finitely many terms vanish). Using the notation of the $q$-Pochhammer symbol, this can be written more efficiently as
\begin{equation} \label{rep_q}
\mathbf{J}_{4_1} (x)  = \sum_{n=0}^\infty |(q;q)_n|^2,
\end{equation}
where $q=e^{2 \pi i x}$. This representation of $\mathbf{J}_{4_1}$ hints at a connection with so-called ``$q$-series'', which play a prominent role in the enumerative combinatorics of partition functions; see for example \cite{berndt}.\\

Coming back to Zagier's problem, clearly the action of $\textup{SL}_2(\mathbb{Z})$ on $\mathbb{Q}$ is generated by the two mappings $x \mapsto x+1$ and $x \mapsto -1/x$. The behavior of $\mathbf{J}_{4_1}$ under the first mapping is trivial by periodicity, but the behavior of $\mathbf{J}_{4_1}$ under the second mapping, i.e.\ the relation between $\mathbf{J}_{4_1}(x)$ and $\mathbf{J}_{4_1}(-1/x)$, is truly fascinating. We already indicate at this point that much of the analysis in the present paper will be based on the theory of continued fractions, which is quite natural since the mapping $x \mapsto 1/x$ plays a central role in that theory. To understand the relation of $\mathbf{J}_{4_1}(x)$ and $\mathbf{J}_{4_1}(-1/x)$, after taking logarithms and switching a sign, Zagier studied the function 
$$
h(x) = \log \frac{\mathbf{J}_{4_1}(x)}{\mathbf{J}_{4_1}(1/x)}, \qquad x \in \mathbb{Q} \backslash \{0\}. 
$$
Since $h(x) = h(-x)$ and $h(x) = -h(1/x)$, it is sufficient to study $h$ on $(0,1)$. Zagier's paper contains several plots of the function $h$, and he writes that the computational data is
\begin{quote}
 ``[\dots] seeming to indicate that the function $h(x)$ is continuous [\dots] at irrational values of $x$.''
\end{quote}

Since $h(x)$ is only defined over rationals, the continuity at irrationals clearly has to be understood with respect to the real topology. In other words, Zagier suggests that $h(x)$ can be extended to a function on $\mathbb{R}$ that is continuous at irrationals. The purpose of the present paper is to prove this conjecture.

\begin{thm} \label{main_thm}
Let $\alpha \in \mathbb{R}$ be irrational. Then the limit $\lim_{x \to \alpha} h(x)$ along rational values of $x$ exists and is finite. In other words: The function $h$ can be extended to a function on $\mathbb{R}$ that is continuous at all irrationals.
\end{thm}

\begin{figure}[ht!]  
\begin{center}
\includegraphics[width=0.8 \linewidth]{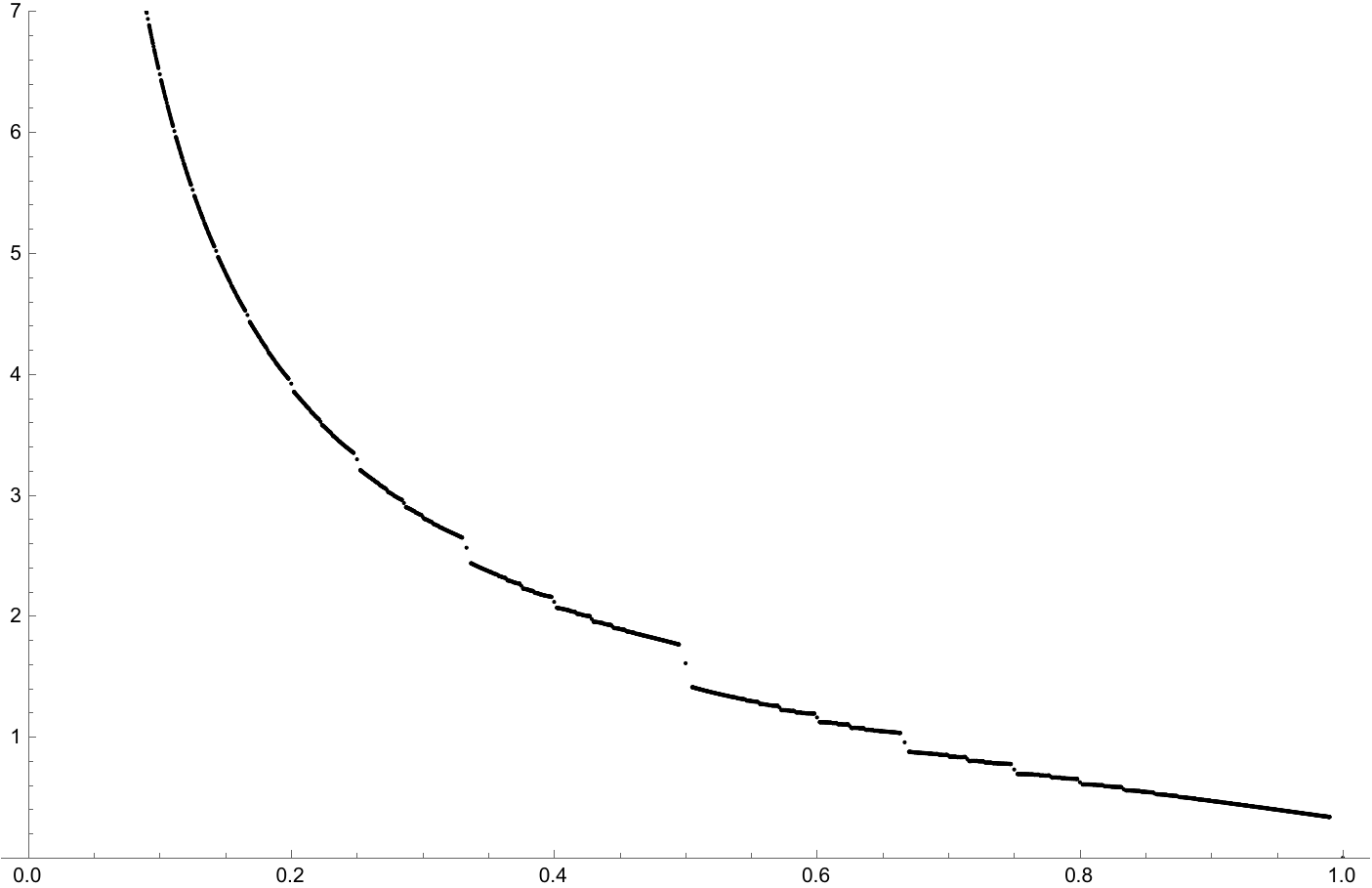}
\end{center}
  \caption{The function $h(x)$, evaluated at all rationals in $(0,1)$ with denominator at most $100$. One can see the relatively big jumps at rationals with small denominators, and the more regular behavior of the function away from such rationals.}  \label{fig:h}
\end{figure}

\begin{figure}[ht!]  
\begin{center}
\includegraphics[width=0.8 \linewidth]{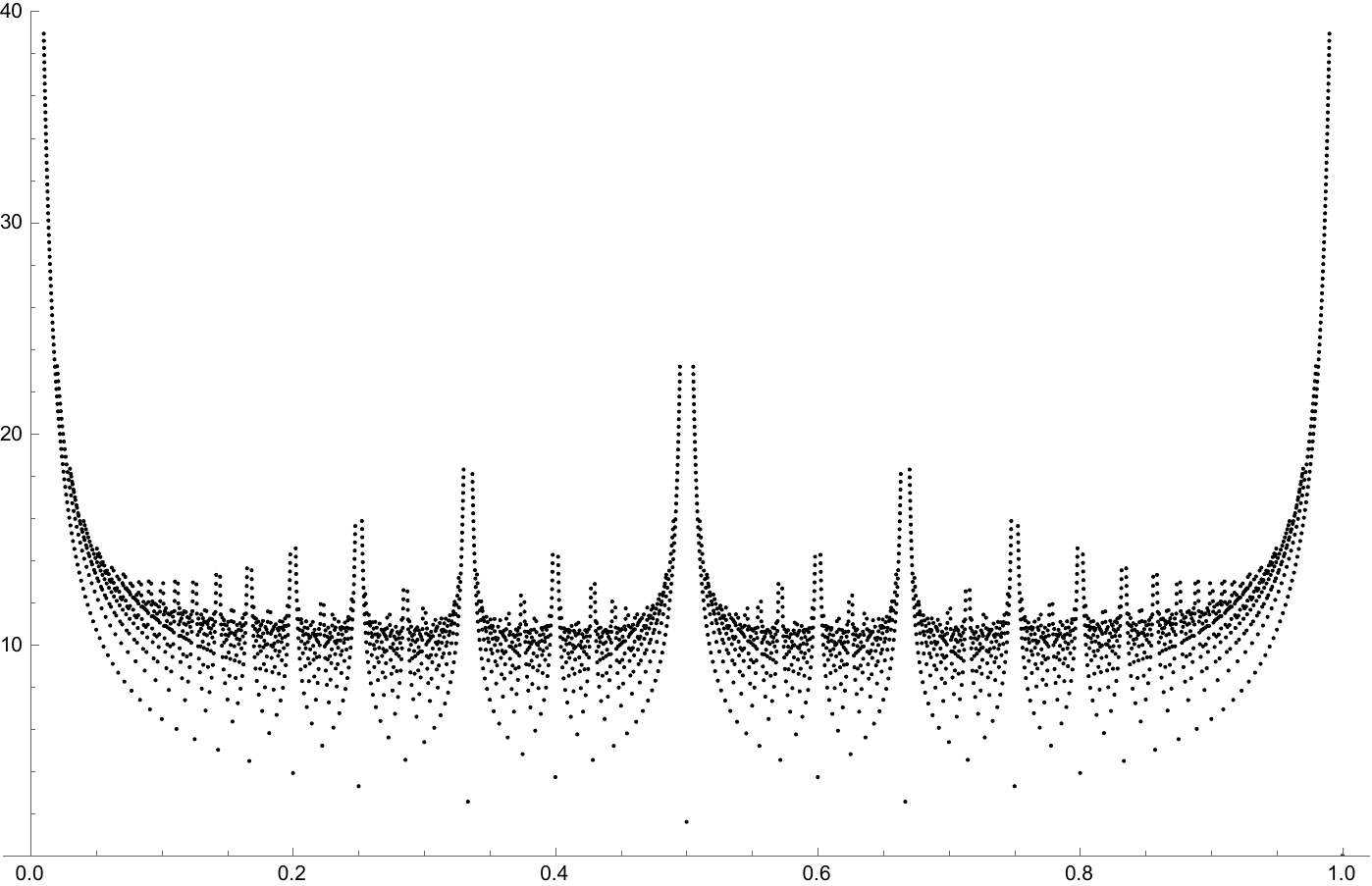}
\end{center}
  \caption{For comparison with the plot of $h(x)=\log \frac{\mathbf{J}_{4_1}(x)}{\mathbf{J}_{4_1}(1/x)}$ in Figure \ref{fig:h}, this is a plot of $\log \mathbf{J}_{4_1}(x)$, again evaluated at all rationals in $(0,1)$ with denominator at most $100$. The plot evidently looks much more irregular than the one in Figure \ref{fig:h}, even if there are indications of a self-similar ``fractal'' structure.}  \label{fig:j}
\end{figure}

A major step towards Theorem \ref{main_thm} was obtained in a work of Aistleitner and Borda \cite{AB3}, where the continuity conjecture was proven for all $\alpha$ satisfying an additional Diophantine property. More precisely, in that paper it was proven that the conclusion of Theorem \ref{main_thm} holds for all irrational numbers $\alpha$ that are not badly approximable, leaving open the continuity of $h$ at badly approximable irrationals. In terms of continued fraction expansions, badly approximable numbers are exactly those that have bounded partial quotients. The assumption of having an \textit{unbounded} sequence of partial quotients in the continued fraction expansion of $\alpha$ played a crucial role for the argument in \cite{AB3}. Quoting from \cite{AB3}: \\ 

``[The main theorem of that paper] leaves the continuity of $h$ at badly approximable irrationals open. It will be seen that our argument crucially relies on the existence of an unbounded subsequence
of partial quotients, so some essential new ideas will be necessary to treat the case of
badly approximable $\alpha$. Some partial results for quadratic irrational $\alpha$ (when the sequence of partial quotients is eventually periodic) are contained in our earlier paper \cite{AB_quantum}. In this case Zagier’s continuity problem might be more tractable than in the general case, due to the additional structure coming from the periodicity of the continued fraction expansion.
The case of general badly approximable $\alpha$ (with no particular structure in the sequence of partial quotients) seems to be even more challenging.''\\

As noted above, in the present paper we prove the conjecture in the fully general case.\\

Broadly speaking, the proof in the present paper is based on methods and on a line of reasoning that are similar to those in \cite{AB3}. The heuristic picture behind the argument of \cite{AB3} is described in detail in \cite[Section 2.3]{AB3}. There, it is also explained why the presence of large partial quotients is crucial for the validity of the argument, since it causes a certain ``independence'' phenomenon that allows us to ``factorize'' the sum in \eqref{q-series} into a product of two sums. In the setup of the present paper, this independence property does not arise automatically from the continued fraction representation of $\alpha$, but instead we distill an ersatz phenomenon out of a statistical analysis of the typical structure of the Ostrowski expansion of positive integers. We will explain the heuristic reasoning behind our proof, and the differences to the argument given in \cite{AB3}, in Section \ref{sec:heu} below, after providing the necessary technical and notational background.\\

\begin{figure}[ht!]  
\begin{center}
\includegraphics[width=0.8 \linewidth]{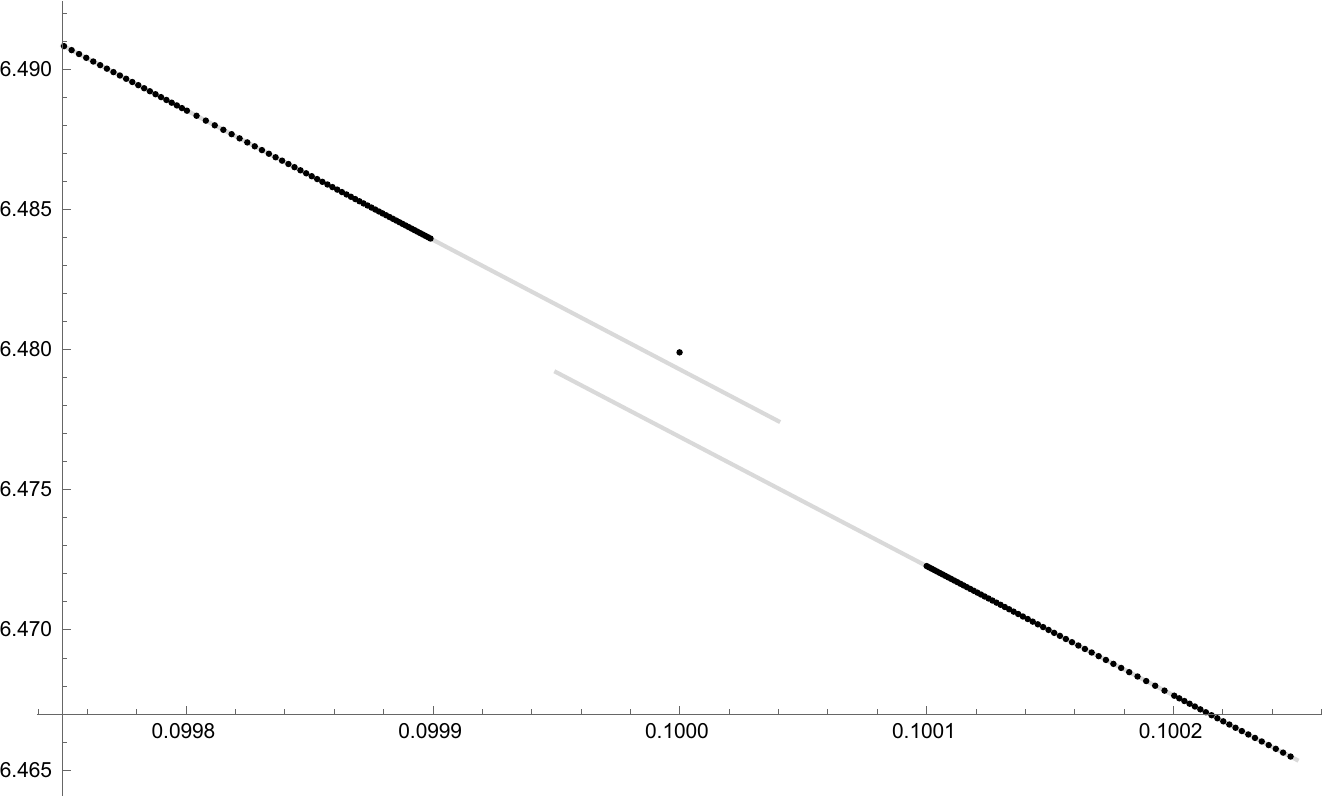}
\end{center}
  \caption{A plot of $h(x)$ in a small neighborhood of $x=1/10$ (function values depicted by black dots; the plot shows the value of $h(x)$ at all rationals $x$  with denominator at most 1000 in the given range), together with the linear regression models for the left and right limits of $h$ at $x$ (solid lines in light gray). The numerical data seem to suggest that $h$ has a (small) upward jump at $x=1/10$, followed by a downward jump, in contrast to the impression of a monotonically decreasing function given by Figure \ref{fig:h}. Note that the existence of left and right limits of $h(x)$ at rational values of $x$ is unproven as of yet.}  \label{fig:jump}
\end{figure}

Before we start with the proofs, we make some concluding remarks. While our paper settles the continuity of $h$ at irrationals, the nature of the jumps of the function at rational arguments remains somewhat mysterious. For example, while the plot of Figure \ref{fig:h} seems to indicate that $h$ is monotonically decreasing in $(0,1)$, with downward jumps to the left and to the right of rationals, the numerical evidence suggests that this is actually not the case; compare Figure \ref{fig:jump}, which shows a plot of $h(x)$ in a small neighborhood of $x=1/10$. The plots also seem to indicate that $h(x)$ has discontinuities at rational values of $x$, but that left and right limits always exist -- all of this remains unproven.

There is no doubt that the proofs in \cite{AB3} and the present paper are designed in an ad-hoc way for the particular case of the colored Jones polynomial of the $4_1$ knot, and are unlikely to allow a generalization to the colored Jones polynomial in the case of other, more complicated, knots. Accordingly, it would be very desirable to have a more conceptual proof of Theorem \ref{main_thm}, which allows a natural generalization to more general knots. On the other hand, the appearance of continued fractions and Ostrowski expansions in our arguments is certainly not purely incidental, and probably the same is true for the involvement of cotangent sums, which play a key role in some of the technical estimates in the proofs, as already in related works in the area of trigonometric products \cite{AB_quantum,other_AB,betdr1,lub}. There exist formulas for $\mathbf{J}_K$ for other hyperbolic knots $K$, which are somewhat similar to those for $\mathbf{J}_{4_1}$ in \eqref{q-series} and \eqref{rep_q}, but more complicated; for example, for the next simplest knot, the $5_2$ (``three-twist knot''), we have (see e.g.\ Formula (2.3) in \cite{kash})
\begin{equation} \label{52_knot}
\mathbf{J}_{5_2}  (x) = \sum_{n=0}^{\infty} \sum_{m=0}^n \frac{(q;q)_n^2}{(\bar{q};\bar{q})_m} q^{-m (n+1)},
\end{equation}
where again $q = e^{2 \pi i x}$ for $x \in \mathbb{Q}$, and $\bar{q}$ denotes the complex conjugate. From a technical perspective, an exceptional property of $\mathbf{J}_K$ in the particular case of $K = 4_1$ is that in this case all summands in \eqref{q-series} and \eqref{rep_q} are positive and real; for other knots $K$, this generally is no longer the case, and there is a high degree of cancellation in the summation formulas such as \eqref{52_knot}, making all calculations extremely delicate; cf.\ \cite{betdr1}. We also note that in Zagier's original paper \cite{ZA} on quantum modular forms, the quotient $\mathbf{J}_K(\gamma x) / \mathbf{J}_K(x)$ was introduced to ``smoothen out'' the rather erratic behavior of $\mathbf{J}_k$ itself, resulting in the function $h$ that (somewhat unsatisfactorily) is still not nicely analytic. Garoufalidis and Zagier \cite{gaza} have then moved on from considering the quotient $\mathbf{J}_K(\gamma x) / \mathbf{J}_K(x)$ to rather upgrading $\mathbf{J}_K$ to a matrix, and reading the quantum modular behavior of $\mathbf{J}_K$ in terms of a matrix product, which leads to a very nice smooth outcome (see also \cite{spirit} for a recent paper which takes this perspective on quantum modularity). It remains open how the results from the present paper align with this matrix perspective on quantum modularity. Finally, we mention that Borda \cite{bence_sawtooth} studied a variant of $\mathbf{J}_{4_1}(x)$ where the q-Pochhammer symbols were replaced by products arising from the periodic sawtooth function. He observed quantum modular behavior of these objects, and formulated  conjectures in the spirit of Zagier's conjecture studied in the present paper. His work indicates that quantum modular behavior might appear in manifold ways in the wider framework of Birkhoff sums for irrational rotations, not necessarily arising from a topological background.

%\begin{equation} \label{formula_h}
%h(x) = \frac{\textup{Vol}(4_1)}{2 \pi x} - \frac{3 \log x}{2} - \frac{\log 3}{4} + o(1)
%\end{equation}
%as $x \to 0^+$ along rationals with bounded numerators; this is in accordance with the numerical data (cf.\ Figure \ref{fig:h}), which in fact seems to suggest that the same holds as $x \to 0^+$ along all rationals. 

\section{Preliminaries} \label{sec:prelim}

As already indicated, Diophantine approximation and the theory of continued fractions will play a key role throughout this paper. We only establish some notation and recall some of the most fundamental facts. For more basic information on continued fractions, we refer to one of the classical texts on the subject, such as those of Khintchine \cite{khin}, Niven \cite{niven}, Schmidt \cite{schmidt} or Rockett and Sz\"{u}sz \cite{rocks}.\\

Throughout the paper, $\alpha$ denotes an irrational real number, and $r$ denotes a rational number.  We write $\alpha=[a_0;a_1,a_2,\ldots ]$ for the (infinite) continued fraction expansion of $\alpha$. The continued fraction expansion is unique, and the positive integers $a_0, a_1, a_2, \dots$ are called the partial quotients of $\alpha$. We write $p_{\ell}/q_{\ell}=[a_0;a_1,\dots, a_{\ell}]$ for the convergents to $\alpha$.  If $r$ is rational, then the continued fraction expansion of $r$ is finite, and we write it as $r = [c_0;c_1, \dots, c_L]$. To make the notation well-defined, we always use the shorter of the two possible continued fraction expansions, namely the one for which $c_L>1$. The number $L$ is called the ``length'' of the continued fraction. The convergents to $r$ are $p_\ell/q_\ell = [c_0;c_1, \dots, c_\ell]$ for $\ell \leq L$, so that $p_L/q_L = r$.\\

We will also need the theory of the Ostrowski numeration system. This is a generalization of the more well-known Zeckendorf numeration system, where integers are represented as sums of Fibonacci numbers under a certain ``digital'' restriction (no two consecutive 1's are allowed). In the Ostrowski system, the denominators of the convergents of some $\alpha$ play the role that the Fibonacci numbers play in the Zeckendorf system (which are the denominators of the convergents in the special case when $\alpha$ is the Golden Mean). Let $\alpha = [a_0;a_1,a_2, \dots]$ be fixed. Then any integer $0 \le N <q_\ell$ has a unique Ostrowski expansion $N=\sum_{i=0}^{\ell-1} b_{i}(N) q_{i}$, where $0 \le b_0(N) < a_1$ and $0 \le b_{i}(N) \le a_{i +1}$ are integers that satisfy the extra rule of $b_{i}(N) =0$ whenever $b_{i +1}(N)=a_{i +2}$. Throughout the paper, we will refer to the coefficients $b_i$ as ``digits'', even if this might be a slight abuse of terminology. Ostrowski numeration is defined analogously for rational $r$ instead of irrational $\alpha$, with the difference that Ostrowski numeration with respect to $\alpha$ is a numeration system on all of $\mathbb{N}$ (by choosing $\ell$ as large as necessary), while Ostrowski numeration with respect to $r= p_L/q_L$ is a numeration system on $\{0, \dots, q_L - 1\}$. \\

Throughout this paper, we will interpret the product on the right-hand side of \eqref{q-series} as a so-called ``Sudler product'', which is a trigonometric product of the form 
\[ P_N (x) := \prod_{n=1}^N |2 \sin (\pi n x)|. \]
With this notation, we have
$$
P_N(x) = \left| (1 - e^{2 \pi i x}) (1 - e^{2 \pi i 2 x}) \cdots (1 - e^{2 \pi i N x}) \right|,
$$
so that for rational $x = p/q$, Equation \eqref{q-series} becomes
\begin{equation} \label{J_as_sum}
\mathbf{J}_{4_1}(p/q) = \sum_{N=0}^{q-1} P_N (p/q)^2.
\end{equation}
Sudler products have a long history going back at least to a paper of Erd\H{o}s and Szekeres \cite{erds}. Among their most interesting aspects are certain self-similarity properties \cite{veme}, which are related to the decomposition in Equation \eqref{pnproductform} below, their relation to cotangent sums, which are known to have a rich arithmetic structure \cite{beco, lub}, and their connection with the spectral theory of almost Mathieu operators \cite{aj,ajm}. A further very interesting connection was developed in Bettin and Drappeau's work on statistics for the distribution of partial quotients of continued fractions, see \cite{betdr2,betdr1}. For a recent survey on Sudler products and generalizations, see \cite{lub_survey}.\\

We will need a shifted form of the Sudler product, which is defined when the number of factors is a convergent denominator of $x$, say the $i$-th denominator $q_i$. Then we define 
$$
P_{q_i}(x,y) := \prod_{n = 1}^{q_i} \left\lvert 2 \sin\Big(\pi\Big(n x + (-1)^i\frac{y}{q_i}\Big)\Big) \right\rvert.
$$
A key technical tool is the decomposition of the full Sudler product $P_N(x)$ into shorter, more controllable, shifted products, according to the Ostrowski decomposition of $N$ with respect to $x$, given by
\begin{equation}\label{pnproductform}
P_N (x ) = \prod_{i =0}^{\ell-1} \prod_{s=0}^{b_{i}(N)-1} P_{q_{i}} \left( x, \varepsilon_{i,s}(N) \right);
\end{equation}
see  \cite[Lemma 2]{AB_quantum} and Proposition \ref{prop_shifted} below. Here, $\ell$ is chosen such that $q_{\ell-1} \leq N < q_\ell$, and 
\begin{equation}\label{epsilondef}
\varepsilon_{i,s}(N) := q_i\left(s \|q_i x\|  + \sum_{j=1}^{\ell-1-i} (-1)^jb_{i+j} \|q_{i+j} x\| \right)
\end{equation}
for $s = 0, \dots, b_i(N)-1$. In the formulas above, and throughout the rest of this paper, we write $\| \cdot \|$ for the distance to the nearest integer. These decomposition formulas hold for all $N$ when $x$ is irrational, and for $N < q$ when $x=p/q$ is rational.

\section{The heuristic picture} \label{sec:heu}

Let $r \in \mathbb{Q} \cap (0,1)$ with continued fraction expansion $[0;c_1,\dots, c_L]$, so that $r = p_L/q_L$. Then the continued fraction expansion of $1/r$ is $[c_1; c_2, c_3, \dots, c_L]$, and since $\mathbf{J}_{4_1}$ is periodic with period 1, what we need to study is the quotient
\begin{equation} \label{this_limit}
h(r) = \log \frac{\mathbf{J}_{4_1}(r)}{\mathbf{J}_{4_1}(r')}
\end{equation}
with $r' = \{1/r\} = [0;c_2, c_3, \dots, c_L]$, where $\{ \cdot \}$ denotes the fractional part. When taking a limit $r \to \alpha$ along rationals, the continued fraction expansion of $r$ ``converges'' to the infinite continued fraction expansion of $\alpha$ (that is, more and more partial quotients at the initial parts of the respective expansions coincide), and the number $r$ in the numerator of \eqref{this_limit} always has the extra partial quotient $c_1$ at the beginning of its continued fraction expansion, in comparison with $r'$ in the denominator.\\

The difficulty of treating the quotient \eqref{this_limit} is that it is a quotient of two sums; recall that according to \eqref{J_as_sum}, we have
\begin{equation} \label{J_as_sum_b}
\mathbf{J}_{4_1} (r) = \sum_{N=0}^{q_L-1} P_N(r)^2,
\end{equation}
where according to \eqref{pnproductform} the product $P_N(r)$ has the factorization
\begin{equation} \label{factor_heur}
P_N(r) = \prod_{i =0}^{\ell-1} \prod_{s=0}^{b_{i}(N)-1} P_{q_{i}} \left( r, \varepsilon_{i,s}(N) \right)
\end{equation}
in terms of the Ostrowski representation of $N$ with respect to the Ostrowski numeration system generated by $r$. Instead of reading the sum in \eqref{J_as_sum_b} as a sum over integers $N < q_L$, we may read it rather as a sum over all possible Ostrowski representations $(b_0, b_1, \dots, b_{L-1})$ of integers $N < q_L$. Therefore we can write
\begin{equation} \label{decompose_this}
\mathbf{J}_{4_1} (r) = \sum_{(b_0, b_1, \dots, b_{L-1})}  \prod_{i =0}^{L-1} \prod_{s=0}^{b_{i}(N)-1} P_{q_{i}} \left( r, \varepsilon_{i,s}(N) \right),
\end{equation}
where the sum ranges over all possible Ostrowski representations $(b_0,b_1, \dots, b_{L-1})$ of integers $N<q_L$. In a similar way, we decompose $\mathbf{J}_{4_1} (r') $ into
$$
\mathbf{J}_{4_1} (r') = \sum_{(b_1, \dots, b_{L-1})}  \prod_{i =0}^{L-1} \prod_{s=0}^{b_{i}(N')-1} P_{q_{i}'} \left( r', \varepsilon_{i,s}(N') \right),
$$
where the sum ranges over all Ostrowski expansions $(b_1, \dots, b_{L-1})$ of integers $N' < q_L'$ (where $q_1', \dots, q_L'$ denote the convergent denominators of $r'$). To see why this can be useful, note that the two sequences of partial quotients of $r$ and $r'$ are very similar, and that accordingly, they generate two closely related systems of Ostrowski numeration (where $r$ has one additional partial quotient $c_1$, and the associated numeration system requires/allows an additional digit $b_0$).\\

Now the key point of the argument in \cite{AB3} is as follows. Assume that it is possible to find an index $k$ between $1$ and $L$ such that we can decompose \eqref{decompose_this} into a product
\begin{eqnarray}
\mathbf{J}_{4_1} (r) & \approx & \left(\sum_{(b_0, b_1, \dots, b_{k-1})}  \prod_{i =0}^{k-1} \prod_{s=0}^{b_{i}(N)-1} P_{q_{i}} \left( r, \varepsilon_{i,s}(N) \right) \right)  \times \nonumber\\
& & \qquad \times \left( \sum_{(b_k, b_{k+1}, \dots, b_{L-1})}  \prod_{i =k}^{L-1} \prod_{s=0}^{b_{i}(N)-1} P_{q_{i}} \left( r, \varepsilon_{i,s}(N) \right) \right) \nonumber\\
& =: & A_k(r) B_k(r), \label{decomp}
\end{eqnarray}
where the first factor depends only on the initial part of an Ostrowski expansion, and the second factor depends only on the tail part of an Ostrowski expansion. Assume that we can similarly decompose 
$$
\mathbf{J}_{4_1} (r')  \approx \left(\sum_{(b_1, \dots, b_{k-1})}  \cdots \right)  \times \left( \sum_{(b_k, b_{k+1}, \dots, b_{L-1})}  \cdots \right) =: A_k'(r') B_k'(r'). 
$$

Assume $k$ to be fixed for the moment. If $r \to \alpha$, then in the products $P_{q_i}(r, \dots)$ appearing in the definition of $A_k(r)$ one can replace $r$ by $\alpha$ with a very small error, since these (finitely many) products depend continuously on $r$. Similarly, one can replace $r'$ by $\alpha' := \{1/\alpha\}$ in the products appearing in the definition of $A_k'(r')$. Thus the quotient $A_k(r)/A_k'(r')$, which depends on $k$ and $r$, can be replaced by a quotient that depends only on $k$ and $\alpha$, and (using the Cauchy convergence criterion) one can show that this quotient converges as $k \to \infty$; this reflects the fact that the influence which the extra partial quotient $c_1$ in $r$ has on the quotient $ \frac{\mathbf{J}_{4_1}(r)}{\mathbf{J}_{4_1}(r')}$ ``stabilizes'' as $r \to \alpha$. For the quotient $B_k(r) / B_k'(r')$, one can see that (in contrast to the formulas for $A_k(r)$ and $A_k'(r')$) the sums for $B_k(r)$ and $B_k'(r')$ both range over the same set of possible Ostrowski digits $(b_k, \dots, b_{L-1})$. Accordingly, one can bijectively map the summands in $B_k(r)$ with the summands of $B_k'(r')$, and show that $B_k(r)/B_k'(r') \to 1$ as $r \to \alpha$. Overall, this proves that $\frac{A_k(r)B_k(r)}{A_k'(r') B_k'(r')}$ converges as $r \to \alpha$, and it is no problem to retain this convergence for $k \to \infty$.\\

The crucial ingredient is to show that a factorization as in \eqref{decomp} is actually possible. There are two difficulties to overcome. Firstly, unlike numeration systems such as the decimal system, whose digits are ``independent'' in an appropriate sense (different decimal digits are stochastically independent with respect to the normalized counting measure on a set such as $\{0, 1, \dots, 10^m -1\}$, for some positive integer $m$), the Ostrowski numeration system has a built-in dependence structure for its digits, which arises as a consequence of the extra rule that $b_{i}(N) =0$ whenever $b_{i +1}(N)=c_{i +2}$ (in stochastic terms, the Ostrowski numeration system does not have independent digits, and instead the digits have a Markov chain structure; see \cite{drst}). This structural dependence of the digit system makes a factorization such as \eqref{decomp} difficult. Secondly, the terms $\varepsilon_{i,s}(N)$ in the shifted products depend on the Ostrowski digits $b_{i+1}, \dots, b_{L-1}$ of $N$; thus all the terms $\varepsilon_{i,s}(N)$ in the second part of the factorization are unproblematic, since they only depend on the digits $(b_k, \dots, b_{L-1})$ covered by that part of the factorization, but the  terms $\varepsilon_{i,s}(N)$ in the first part of the factorization are problematic, since they also depend on the digits that are only supposed to enter the second part of the factorization. These are two genuinely different problems, but both of them could be settled in \cite{AB3} thanks to the assumption of the existence of arbitrarily large partial quotients in the continued fraction expansion of $\alpha$. We refer to \cite[Section 2.3]{AB3} for a more detailed exposition, but roughly speaking, both problems can be solved if the factorization \eqref{decomp} is carried out at an index $k$ such that the following partial quotient $c_{k+1}$ is ``very large''.\\

In the setup of the present paper, we are not provided with the existence of such ``very large'' partial quotients. Accordingly, we must find a different solution for the two problems described in the previous paragraph. We note that the first problem (dependence of the digits within the Ostrowski numeration system) only arises when a digit $b_{i+1}$ attains its maximal potential value, forcing the preceding digit $b_i$ to be $0$; assuring that $b_{i+1}=0$ would break this dependence between the digits with index smaller than $i+1$ and those with index larger than $i+1$. This idea of breaking dependencies in the Ostrowski numeration system for typical Ostrowski expansions was already exploited in a recent paper of the second-named author in order to establish equidistribution in $\mathbb{Z}_d$ in certain Bohr sets arising in Diophantine approximation \cite{bad_rough}.\\

Concerning the second problem, note that as the formula for $\varepsilon_{i,s}(N)$ shows, these numbers depend on the Ostrowski digits $b_{i+1}, b_{i+2} , \dots$ in a complicated way, but in such a way that $b_{i+1}$ typically contributes most to $\varepsilon_{i,s}(N)$, while $b_{i+2}$ contributes less, $b_{i+3}$ contributes even less, and so on. If we could ensure that $b_{i+1}= b_{i+2} = \dots = b_{i+m} = 0$ for some (sufficiently large) $m$, then this would make the contribution of the sum over $j$ in formula \eqref{epsilondef} very small and essentially yield $\varepsilon_{i,s}(N) \approx s q_i \|q_i \alpha\|$, thereby resolving the second problem towards a factorization as in \eqref{decomp}. Accordingly, both of our problems can be settled if we can assure that there is a long run of consecutive zeros in the Ostrowski representation. Now, in the setup of the present paper, it is indeed true that for a ``typical'' integer $N < q_L$, we can expect a long run of consecutive zeros in its Ostrowski representation -- here we crucially use the fact that by assumption the partial quotients are bounded, so that each Ostrowski digit only has a finite, uniformly bounded, number of possible values. This ingredient is in the spirit of the Erd\H{o}s--R\'enyi ``pure heads'' theorem, which asserts that when tossing a coin $u$ times, one should expect to see a run of roughly $\log u$ many consecutive heads (see \cite{erdren,erdrev}). In this way, we are able to simultaneously break the dependence structure arising from the Ostrowski numeration system on the one hand, and from the influence of the $\varepsilon_{i,s}(N)$ terms on the other hand. Note, however, that in the situation of \cite{AB3} as described in the previous paragraph, the factorization \eqref{decomp} was carried out at a certain (fixed) index $k$ for which $c_{k+1}$ is ``very large''. In contrast, now we aim at a factorization which is based upon the existence of long runs of zeros in the Ostrowski representation of $N$, but while statistical reasoning ensures that such a long run of zeros exists for most integers $N$, we clearly cannot expect that this long run of zeros always occurs at the \textit{same} location within the digital representation $(b_0, b_1,\dots, b_{L-1})$ of $N$. Accordingly, instead of being able to factorize $\mathbf{J}_{4_1}(r)$ at a fixed index $k$ as in \eqref{decomp} and \cite{AB3}, in the present paper we will apply a factorization along a ``running index'', which accounts for the different possibilities of the location of a long run of zeros in the Ostrowski expansion of $N$.\\

A final remark on the proof. In view of the partial solution provided by \cite{AB3}, throughout the present paper we may assume that $\alpha$ is badly approximable (and thus has bounded partial quotients in its continued fraction expansion). We are interested in the behavior of $h(r)$, as $r$ approaches $\alpha$. As $r \to \alpha$, more and more partial quotients at the initial segment of the continued fraction expansion of $r$ coincide with those of $\alpha$, and thus are also bounded. However, $r \to \alpha$ emphatically does \emph{not} imply that \emph{all} partial quotients of $r$ can be assumed to be bounded -- on the contrary, we must make allowance for the possibility that some later partial quotients of $r$ could be extremely large. Accordingly, throughout the argument the continued fractions / Ostrowski expansions will be split into two segments: an initial part, where the boundedness of the partial quotients of $\alpha$ carries over to the partial quotients of $r$ (and will be crucially used), and a tail part where we have to work in fully general circumstances, without any control of the potential size of the partial quotients of $r$. We are in the fortunate situation that many (highly non-trivial) estimates for the tail part can be adopted directly from \cite{AB3}, since there as well as in the present paper, no assumptions on Diophantine properties related to the tail part can be made.

\section{Proof of Theorem \ref{main_thm}} \label{sec:proof}

We will first introduce the general machinery for the proof of Theorem  \ref{main_thm}, and formulate several auxiliary lemmas. In Section \ref{sub:assuming} we will give the proof of Theorem \ref{main_thm}, assuming the validity of these lemmas. Afterwards, in Section \ref{sec_techn}, we will give the proofs of the lemmas.

\subsection{Admissible tuples} \label{sec:ad}

Let $\alpha = [0;a_1,a_2,\ldots]$ be a badly approximable irrational number, which will be understood to remain fixed throughout the rest of the paper. Let $M = M(\alpha) := \max_{i \in \mathbb{N}}a_i < \infty$. As mentioned above, the Ostrowski expansion of a non-negative integer $N$ is the representation
\begin{equation*}
N = \sum_{i = 0}^\ell b_{i}(N)q_{i} \quad \text{ where } 0 \ls b_0 < a_1, \;\;  0 \ls b_{i} \ls a_{i+1} \text { for } i \gs 1,
\end{equation*}
with the extra rule that
$b_{i-1} = 0$ whenever $b_{i} = a_{i+1}$. 
This representation is unique if the leading digit $b_\ell$ is assumed to be non-zero. We say that a tuple $(b_0,\ldots,b_{K-1})$ is \textit{admissible} (with respect to $\alpha$, which is omitted if clear from the context) if 
\[0 \ls b_0 < a_1, \;\;  0 \ls b_{i} \ls a_{i+1} \quad \text { for } 1 \leq i \leq K-1,\]
and if $b_{i-1} = 0$ whenever $b_{i} = a_{i+1}$; in other words, admissible tuples are those that specify possible Ostrowski expansions of an integer (of given length, and with respect to $\alpha$). We write $\mathcal{A}_K = \mathcal{A}_K(\alpha)$ for the set of all admissible $K$-tuples.
In that way, we can define the bijection
\[\begin{split}\psi_K = \psi_{K,\alpha}: \{0,1,\ldots,q_{K}-1\} &\to \mathcal{A}_K\\
N  &\mapsto (b_0(N),\ldots,b_{K-1}(N))
\end{split}\]
where
$N = \sum_{i = 0}^{K-1} b_{i}(N)q_{i}$ in Ostrowski representation. We extend this to 
all natural numbers: We write $\tilde{\mathcal{A}_K}$ for the subset of $\mathcal{A}_K$ such that $b_{K-1} \neq 0$. In that way, ${\psi_K}$ maps the set $\{q_{K-1},\ldots,q_{K}-1\}$ bijectively to $\tilde{\mathcal{A}_K}$. (For completeness, we also need to define $\tilde{\mathcal{A}_0}$ to be the set of the 1-tuple $(0)$, so that $\{0\}$ is mapped to  $\tilde{\mathcal{A}_0}$ and the integer $N=0$ also is correctly handled; this is a special case, since it corresponds to the only Ostrowski expansion with a leading zero). This allows us to define the bijection
\begin{eqnarray}\psi = \psi_{\alpha}: \mathbb{N} &\to & \mathcal{A} := \dot{\bigcup_{K \in \mathbb{N}}}\tilde{\mathcal{A}_K}, \nonumber\\
N &\mapsto & (b_0(N),\ldots,b_{K-1}(N)), \quad b_{K-1}(N) \neq 0, \label{is_mapped}
\end{eqnarray}
(again with the special case $0 \mapsto (0)$), where the suitable value of $K$ in \eqref{is_mapped} depends on $N$.\\

In a very similar way, we can define admissible tuples and a bijection between integers and Ostrowski expansions with respect to a rational number $r$ (instead of an irrational $\alpha$, as in the previous section). Let $r \in \mathbb{Q}$ be given with finite continued fraction expansion $r = [0;c_1,\ldots,c_L]$, so that $r = p_L/q_L$. We can define sets $\mathcal{A}_K$ of admissible $K$-tuples with respect to $r$ for all $K$ up to $L-1$ analogously to the definitions in the irrational case. We can also define bijections $\psi_K = \psi_{K,r}:~ \{0, 1, \dots, q_K - 1\} \to \mathcal{A}_K$, and sets $\tilde{\mathcal{A}_K} \subset \mathcal{A}_K$, analogous to the above, for all $K \leq L-1$. In the rational case, we do not consider tuples $(b_0, \dots, b_{K-1})$ whose length $K$ exceeds the length $L$ of the continued fraction expansion of $r$. Accordingly, in the rational case we can construct a bijective function
\[\begin{split}\psi = \psi_r: \{0, \dots, q_L - 1\} &\to \mathcal{A} := \dot{\bigcup_{K \leq L}}\tilde{\mathcal{A}_K}, \\
N &\mapsto (b_0(N),\ldots,b_{K-1}(N)), \quad b_{K-1}(N) \neq 0
\end{split}\]
(and $0 \mapsto (0)$), where again the suitable value of $K \leq L$ depends on $N$.\\\\ 

Finally, we define 
\begin{equation} \label{J_def}
J(r) := \sum_{N = 0}^{q_L-1} P_N(r)^2, \qquad \text{where} \qquad P_N(x) := \prod_{n = 1}^N 2|\sin(\pi x)|.
\end{equation}

\subsection{Finding a run of consecutive zeros} \label{sec:findrun}
Let $r$ be a rational number from the interval $(0,1)$. Let $r = [0;c_1,\ldots,c_L]$ be its continued fraction expansion, of (finite) length $L$. We fix $k < L$ and set $m := \lfloor \log \log k\rfloor, t := \lfloor\tfrac{k}{m}\rfloor$. At the end of the proof we will have $L \to \infty$ (as a consequence of $r \to \alpha$) and choose $k$ ``large'', so for simplicity of writing we can assume throughout the paper that $k \geq 20$, say. As a consequence, the integer $m$ from above is well-defined and positive. We emphasize that the Ostrowski expansions appearing in Sections \ref{sec:findrun}--\ref{sub:assuming}  are understood to be taken with respect to $r$, not with respect to $\alpha$. Similarly, throughout these sections $p_\ell/q_\ell$ are convergents to $r$, not convergents to $\alpha$.\\

For $t/2 \leq j < t$, we define
\[\mathcal{F}_{j,L}= \mathcal{F}_{j,L}(k,r) := \left\{0 \leq N < q_{L}-1: b_{i}(N) = 0\quad  \forall i \in \{jm + \ell, 0 \leq \ell < m\} \right\}, \]
and write
\[\mathcal{G}_{j,L} := \mathcal{F}_{j,L}\setminus \bigcup_{t/2 \leq i < j}\mathcal{F}_{i,L}.\]
We say an integer $N < q_L$ is \textit{good} if \[N \in \mathcal{G}_L := \bigcup_{t/2 \leq j < t}\mathcal{F}_{j,L} = \dot{\bigcup_{t/2 \leq j < t}}\mathcal{G}_{j,L},\]
and $N < q_L$ is \textit{evil} if $N \notin \mathcal{G}_L$.
We denote by $\mathcal{E}_L = \mathcal{E}_L(k,\alpha)$ the set of evil numbers up to $q_L$. Roughly speaking, the good integers $N$ are those whose Ostrowski expansion contains a long run of consecutive zeros; note how being contained in $\mathcal{F}_{j,L}$ means that there is a long run of zeros starting at a location specified by the index $j$, while being in $\mathcal{G}_{j,L}$ means that the \emph{first} such long run of zeros starts at this location. We will show in the sequel that being good is a generic property, in the sense that ``most'' integers are good, and that (crucially) it is essentially the contribution of only the good integers which determines the size of $J(r)$ in \eqref{J_def}. Note that in our construction we are only trying to find a long run of zeros among the first $\approx tm \approx k$ Ostrowski digits of $N$, not among all $L$ digits (where later $k$ will be assumed to be large but fixed, while $L \to \infty$). This is because knowing that $\alpha$ is badly approximable and that $r$ is close to $\alpha$ provides us with a bound for the size of the partial quotients of $r$ with small index, but with no control over the size of partial quotients with large index. More precisely, the choice of $t$ and $m$, and the construction of the sets $\mathcal{F}_{j,L}$ and $\mathcal{G}_{j,L}$ above, is made in such a way that the largest index of an Ostrowski digit that is relevant for any $\mathcal{F}_{j,L}$ (resp.\ $\mathcal{G}_{j,L}$) is the last digit relevant for the case $j = t-1$, namely the Ostrowski digit with index $(t-1) m + m -1 = tm -1$. By our choice of $t$ and $m$ we have $t m -1 < k$. Thus if we can guarantee that the boundedness of \emph{all} partial quotients of $\alpha$ carries over to the \emph{initial} $k$ partial quotients of $r$, then we can guarantee that any digit $b_i$ whose index is within our ``finding a long run of zeros'' region can only take a bounded number of possible values, a fact which will be crucially used in the proof of (for example) Lemma \ref{lem_evil_good_1} and Corollary \ref{cor_evil_good_2} below. \\

We now define maps from $\mathcal{E}_L$ into $\mathcal{G}_{j,L}$ by replacing in the Ostrowski expansion of $N$ the digits $b_i, i \in \{jm + \ell, 0 \leq \ell < m\}$ by $0$'s: More precisely, we define the maps
$\pi_j: \mathcal{E}_L \to \mathcal{G}_{j,L}$ via 
$\pi_j := \psi^{-1} \circ p_j \circ \psi$, where
\[\begin{split} p_j: \mathcal{A}_L &\to \mathcal{A}_L\\
(b_0,\ldots,b_{L-1}) &\mapsto (b_0,\ldots,b_{jm -1},0,\ldots,0,b_{m(j+1)},\ldots,b_{L-1}).\end{split}\]
Here it is important to note that $(b_0,\ldots,b_{jm -1},0,\ldots,0,b_{m(j+1)},\ldots,b_{L-1})$ indeed is an admissible tuple, since the Ostrowski rule never forbids replacing non-zero digits with zeros. Further, it is important to note that $\pi_j$ is indeed mapping $\mathcal{E}_L$ to $\mathcal{G}_{j,L}$ (and not only to $\mathcal{F}_{j,L}$): Assuming the converse, there must exist an $i$ with $t/2 \leq i < j$ such that $(b_{im}(\pi_j(N),\ldots,b_{i(m+1) -1}(\pi_j(N)) = (0,\ldots,0)$. Since $p_j$ keeps those digits fixed, this implies 
$(b_{im}(N),\ldots,b_{i(m+1) -1}(N)) = (0,\ldots,0)$,  which yields 
$N \in \mathcal{F}_{i,L}$, a contradiction to $N \in \mathcal{E}_L$.\\

We note that $\pi_j$ is not injective, but it will turn out that only a small number of evil elements can be mapped onto the same element in $\mathcal{G}_{j,L}$. This will be used in the proof of the following statement later on.

\begin{lem}   \label{lem_evil_good_1}  %[Difficulty $8/10$]
For given $r=[0;c_1, \dots, c_L]$, let $k,m,t$ be defined as above, and let $N \in \mathcal{E}_L$. Assume that $\max_{i \leq k+2} c_i \leq M$. Then for every $j$ with $t/2 \leq j < t$ we have
    \[P_N(r) \ll ({\mathcal{O}_M (1)})^{m} P_{\pi_j(N)}(r),\]
    where the implied constants only depend on $M$.
\end{lem}

Lemma \ref{lem_evil_good_1} asserts that whenever $N$ is an evil integer, then for every $t/2 \leq j < t$, we can also find a suitable good integer $\pi_j(N) \in \mathcal{G}_{j,L}$ whose Sudler product $P_{\pi_j(N)}$ is of roughly similar size as $P_N$. Here it is important that we have a whole range for the choice of $j$, so to one specific evil integer $N$ we find not one, but many different good integers contributing to $J(r)$. This is used in Corollary \ref{cor_evil_good_2} below to show that the main contribution to $J(r)$ comes from the good, and not from the evil, integers. 

\begin{cor}\label{cor_evil_good_2}    %[Difficulty $3/10$]
Let $\varepsilon > 0$ be fixed. There exists $K_0=K_0(\varepsilon)$ such that whenever $k > K_0$, we have the following. Assume that $r = [0;c_1,\ldots,c_L]$ satisfies $\max_{i \leq k+2} c_i \leq M$, and let $m$ and $t$ be defined as above. Then we have
\begin{equation} \label{evilgood}
\sum_{N \in \mathcal{E}_L}P_N^2(r) < \varepsilon J(r).
\end{equation}
\end{cor}

\begin{proof}[Proof of Corollary \ref{cor_evil_good_2} assuming Lemma \ref{lem_evil_good_1}]
Applying Lemma \ref{lem_evil_good_1}, we have

\[t \sum_{N \in \mathcal{E}_L} P_N^2(r) \ll \sum_{t/2 \leq j < t} \sum_{N \in \mathcal{E}_L}  P_{N}^2(r) \ll {({O_M(1)})}^{m} \sum_{t/2 \leq j < t} \sum_{N \in \mathcal{E}_L}  P_{\pi_j(N)}^2(r).\]
    We observe that $|\pi_j^{-1}(N)| \leq (M+1)^m$, since for any particular admissible tuple $(b_0,\ldots,b_{jm -1},0,\ldots,0,b_{(j+1)m},\ldots,b_{L-1})$, there are at most $M+1$ possibilities to replace any specific zero digit in $({jm},\ldots,{jm + m-1})$ by some other digit $d_i \in \{0,\ldots,c_{i+1}\}$ (without violating admissibility), since $c_{i+1} \leq M$ for all relevant $i$ by assumption. Here we crucially used (as explained at the beginning of this section) that by our choice of $t$ and $m$, the maximal possible index $i$ of a digit $b_i$ that is changed by some $\pi_j$ is of size $i = (t-1) m + m - 1 = tm-1$ (corresponding to the case $j = t-1$), and we have $tm - 1 < k$ so that indeed $d_i \leq c_{i+1} \leq M$ for all indices $i$ that are relevant to this argument. Using Lemma \ref{lem_evil_good_1}, noting that $\mathcal{G}_{j,L} \cap \mathcal{G}_{j',L} = \emptyset$ for $j \neq j'$, we have

    \[\begin{split} t\sum_{N \in \mathcal{E}_L}  P_{N}^2(r) &\ll {({O_M(1)})}^m 
    \sum_{t/2 \leq j < t} \sum_{N \in \mathcal{E}_L}P_{\pi_j(N)}^2(r) 
    \\&\ll {({O_M(1)})}^{m} \sum_{t/2 \leq j < t} \sum_{N \in \mathcal{G}_{j,L}}P_{N}^2(r)|\pi_j^{-1}(N)|  \\&\ll {({O_M(1)})}^{m}\sum_{N < q_L}P_N^2(r),\end{split}\]
    so that 
    $$
    \sum_{N \in \mathcal{E}_L}  P_{N}^2(r)  \ll \frac{{({O_M(1)})}^{m}}{t} \sum_{N < q_L}P_N^2(r).
    $$
Since $t \sim \frac{k}{\log \log k}$ and $m \sim \log \log k$, this finishes the proof.
\end{proof}

Roughly speaking, what Corollary \ref{cor_evil_good_2} asserts is the following: There are some evil integers $N$, but to each such evil integer we can associate a large number (namely: order $t$ many) good integers. The value of $P_N(r)$ for the evil $N$ might slightly exceed the value of the Sudler product of the associated good integers (Lemma \ref{lem_evil_good_1}), but this is compensated by the fact that to each evil $N$ we can associate a large number of good integers. On the other hand, each good integer is associated only with a limited number of evil integers, as a consequence of the boundedness of the Ostrowski digits (which comes from the boundedness of the partial quotients). Accordingly, the main contribution to $J(r)$ comes from the good, and not from the evil integers, as witnessed by \eqref{evilgood}. Note that the boundedness of the partial quotients has to be used twice: once essentially to compare the cardinality of $\mathcal{E}_L$ with the cardinality of $\mathcal{G}_L$ (very easily, in the ``proof of Corollary \ref{cor_evil_good_2} assuming Lemma \ref{lem_evil_good_1}'', where we also exploit the fact that we have many options for the index $j$ that localizes a long run of zeros), and once (in a much more fundamental way) in the proof of Lemma \ref{lem_evil_good_1} to compare the size of $P_N(r)$ with that of $P_{\pi_j(N)}(r)$.

\subsection{The splitting process for good $N$}
Let $t/2 \leq j < t$ fixed. We define $\mathcal{G}_{j,L}^{(1)}$ as $\mathcal{G}_{j,L} \cap \{1,\ldots,q_{jm}-1\}$, and  note the crucial observation that $\mathcal{G}_{j}^{(1)} := \mathcal{G}_{j,L}^{(1)}$ does not actually depend on $L$ (which is determined by $r$), but only depends on our choice of $k$, since by construction $jm \leq k < L$. Further, we define $\mathcal{G}_{j,L}^{(2)}$ (which now does indeed depend on $L$) by 
\[\mathcal{G}_{j,L}^{(2)} := \{N_2 < q_L: b_i(N_2) = 0 \quad \forall i < (j+1) m\}.\] Roughly speaking, since $\mathcal{G}_{j,L}$ contains integers whose Ostrowski representation has a long run of zeros starting at index $jm$, in $\mathcal{G}_{j}^{(1)}$ we encode those digits that come before the run of zeros (i.e.\ digits with small index), and in $\mathcal{G}_{j,L}^{(2)}$ we encode the digits that come afterwards. Note that since the digits with small index and those with large index are separated by zeros, indeed every initial segment can be combined with every tail segment, since the ``extra rule'' of the Ostrowski numeration system does not apply. In mathematical terms, there is a bijection
\begin{equation}\label{bij_theta}\begin{split}\theta_j: \mathcal{G}_{j}^{(1)} \times \mathcal{G}_{j,L}^{(2)} &\to \mathcal{G}_{j,L}\\
(N_1,N_2) &\mapsto N_1 + N_2.
\end{split}\end{equation}
By inverting this bijection, to each $N \in \mathcal{G}_{j,L}$ we can assign unique numbers $N_1 \in \mathcal{G}_{j}^{(1)}$ and $N_2 \in \mathcal{G}_{j,L}^{(2)}$ such that $N = N_1 + N_2$. With this notation at hand, we will  prove that the Sudler products $P_N$ for $N \in \mathcal{G}_{j,L}$ decompose approximately into $P_{N_1}\cdot P_{N_2}$ (individually, before taking a summation over $N$), and that accordingly the sum $\sum_N P_N$ decomposes into a product $\sum_{N_1} P_{N_1}(r) \times \sum_{N_2} P_{N_2}(r)$. We emphasize once more that the existence of a run of zeros in the Ostrowski expansion of $N$ is used twice: Firstly, by breaking the dependence in the Ostrowski numeration system (coming from the ``extra rule'' of this numeration system) to guarantee that the index set $\mathcal{G}_{j,L}$ indeed decomposes into a product $\mathcal{G}_{j}^{(1)} \times \mathcal{G}_{j,L}^{(2)}$, and secondly, to guarantee that $P_N(r) \approx P_{N_1}(r)\cdot P_{N_2}(r)$, essentially by breaking the dependence structure of the Sudler products, which is encoded in the ``shifts'' $\varepsilon_i$ in the factorization formula \eqref{factor_heur}. These are two different effects. The first one needs only one $0$ digit, the second one a sufficiently long run of zeros. The first one is tied to the way how $J(r)$ arises as a sum over integers, interpreted as a sum over configurations of admissible digits, while the second one is a ``pointwise'' effect which holds for particular individual values of $N$.

\begin{lem}  \label{lem_decomp}   %[Difficulty 7/10]
Let $r = [0;c_1,\ldots,c_L] \in \mathbb{Q}$ such that $\max_{1 \leq i \leq k+2} c_i \leq M$ for a fixed $k \leq L-2$. For all $k$ there exists $\eta_k > 0$ with $\eta_k \to 0$ as $k \to \infty$, such that for all $j$ in the range $t/2 \leq j < t$, we have
    \[\frac{\sum_{N \in \mathcal{G}_{j,L}}P_N^2(r)}{\sum_{N_1 \in \mathcal{G}_{j}^{(1)}}P_{N_1}^2(r) \cdot \sum_{N_2 \in \mathcal{G}_{j,L}^{(2)}}P_{N_2}^2(r)} \in (1 - \eta_k, 1 + \eta_k).\]
  The numbers $\eta_k$ do not depend on $j$, and are uniform among all $r$ for which the first $k$ partial quotients coincide.
\end{lem}

\begin{lem} \label{lem_transf_princ}   %[``gratis'']
For $D \in \mathbb{N}$, let $I_D$ denote the interval around the irrational $\alpha \in (0,1)$ which consists of all numbers whose (finite or infinite) continued fraction expansion also starts with the segment $[0;a_1,\dots,a_D]$. Then for every $k \in \mathbb{N}$ and every $\varepsilon > 0$, there exists $D = D(k,\varepsilon,\alpha) \in \mathbb{N}$ such that
\[\max_{N < q_k}\sup_{\beta \in I_D}\left\lvert\frac{P_N(\beta)}{P_N(\alpha)} - 1\right\rvert < \varepsilon.\]
Here $q_k$ denotes the $k$-th convergent denominator of $\alpha$. The same holds true if $\alpha$ is not an irrational, but a rational with denominator greater than $q_k$. 
\end{lem}

\begin{proof} This is just the fact that $P_N(\beta)$ depends on $\beta$ in a continuous way, and that $P_N(\alpha)$ is non-zero for irrational $\alpha$ (resp.\ for rational $\alpha$ with denominator greater than $N$), together with the fact that the length of $I_D$ goes to zero as $D \to \infty$.

% Altes quantitatives Argument: 
%The above statement is clearly equivalent to demanding (by altering $\varepsilon$)
%\[\sup_{N < q_k}\sup_{\alpha \in I_D}\left\lvert\log {P_N(r)} - \log{P_N(\alpha)}\right\rvert < \varepsilon.\]
%Writing $f(x) = \log(2\sin \pi x)$, we have by the intermediate value theorem (and using the fact that $\alpha$ is badly approximable)
%\[\left\lvert\log {P_N(r)} - \log{P_N(\alpha)}\right\rvert
%\leq \sum_{n \leq q_k} \lvert f(nr) - f(n\alpha) \rvert
%\leq \left\lvert f'(q_k\alpha + O(q_D^{-2})) \right\rvert \sum_{n \leq q_k} \lvert n(\alpha -r)\rvert
%\ll_{\alpha} \frac{q_k^4}{q_D^{2}}.\]
%Choosing now $D$ sufficiently large in terms of $k$ immediately proves the claim.

\end{proof}

 We define $r' := \{1/r\} = [0;c_2,\ldots,c_L]$, and set $p_{i}'/q_{i}' = [0;c_2,\ldots,c_i]$ for $i \geq 1$ (note that this is the $(i-1)$-th convergent to $r'$, with $p_{1}'/q_{1}' := 0/1$). In particular, $r' = p_L'/q_L'$. We also define the function $S = {\pi^{-1}_{r'}} \circ s \circ {\pi_{r}}$, 
where $s((b_0,b_1,\ldots,b_i)) := (b_1,\ldots,b_i)$. Essentially, $S$ arises from a shift on the Ostrowski digits, and maps $\{0, \dots, q_L -1\}$ to $\{0, \dots, q_L' -1\}$, in such a way that $N = \sum_{i = 0}^{L-1}b_{i}(N)q_{i}$ is mapped to $N' := S(N) = \sum_{i = 1}^{L-1}b_{i}(N)q_{i}'$. The numbers $N$ and $N'$ are related by the fact that the Ostrowski expansion of $N$ with respect to $r$ is the same as the Ostrowski expansion of $N'$ with respect to $r'$, except for the extra digit $b_0(N)$ of $N$. We will use this to relate the value of the Sudler product $P_N(r)$ to that of $P_{N'}(r')$, which is plausible in view of the factorization \eqref{factor_heur}.

\begin{lem}    \label{lem_tail_conv}   %[Difficulty 2/10]
     For all $k$ there exists $\eta_k > 0$ such that $\eta_k \to 0$ as $k \to \infty$, and such that for all $r \in \mathbb{Q}$ with $L > k+1$ and $\max_{1 \leq i \leq k+2}c_i \leq M$, we have
    \[\frac{\sum_{t/2 \leq j < t}\sum_{N_2 \in \mathcal{G}_{j,L}^{(2)}(k,r)}P_{N_2}(r)}{
    \sum_{t/2 \leq j  < t}\sum_{N_2 \in S(\mathcal{G}_{j,L}^{(2)}(k,r))}
    P_{N_2}(r')} \in (1 - \eta_k, 1+\eta_k).\]
\end{lem}

Finally, we need the following analogues of Corollary \ref{cor_evil_good_2} and Lemma \ref{lem_decomp} for $r'$ instead of $r$. The proofs are the same as the ones given above, apart from the fact that the run of consecutive $0$ for ``good numbers'' now starts exactly one position earlier (since the digit $b_0$ disappeared when switching from $r$ to $r'$).

\begin{cor}[Corollary \ref{cor_evil_good_2} for $r'$]\label{cor_r'}
Let $\varepsilon > 0$ be fixed. There exists $K_0$ such that whenever $k > K_0$, we have the following.
Assume that $r = [0;c_1,\ldots,c_L]$ and $r' = [0;c_2,\ldots,c_L]$ are such that $\max_{i \leq k+2} c_i \leq M$.  Let
$S$ be the mapping from above. Then we have
    \[\sum_{N \in S(\mathcal{E}_L)}P_{N}^2(r') < \varepsilon J(r').\]   
\end{cor}

\begin{lem}[Lemma \ref{lem_decomp} for $r'$]\label{split_r'}
Let $r \in \mathbb{Q}$. For all $k$ there exists $\eta_k > 0$ such that $\eta_k \to 0$ as $k \to \infty$, and such that the following holds. If $L > k$, then for all $j$ in the range $t/2 \leq j < t$ we have
    \[\frac{\sum_{N \in S(\mathcal{G}_{j,L})}P_N^2(r')}{\sum_{N_1 \in S(\mathcal{G}_{j}^{(1)})}P_{N_1}^2(r') \cdot \sum_{N_2 \in S(\mathcal{G}_{j,L}^{(2)})}P_{N_2}^2(r')} \in (1 - \eta_k, 1 + \eta_k).\]
  The numbers $\eta_k$ do not depend on $j$, and are uniform among all $r$ for which the first $k$ partial quotients coincide.
\end{lem}

\subsection{Proof of Theorem \ref{main_thm} assuming technical Lemmas} \label{sub:assuming}

In this section, we assume that Lemmas \ref{lem_evil_good_1}, \ref{lem_decomp}, \ref{lem_tail_conv} and  \ref{split_r'} (and thus also Corollaries \ref{cor_evil_good_2} and \ref{cor_r'}) are all true. We will show how they imply Theorem \ref{main_thm}. The auxiliary results will then be proven in Section \ref{sec_techn}.\\

As in earlier parts of this section, we assume that $\alpha$ is badly approximable and consequently, $M := \max_{i \in \mathbb{N}} a_i(\alpha)$ is finite. We consider a rational $r \in (0,1)$ and $r ' = \{1 / r\}$, and study 
\[h(r) = \log \left( \frac{\sum_{N=0}^{q_L-1} P_N(r)^2}{\sum_{N=0}^{q_L'-1} P_N(r')^2}\right) \qquad \text{as $r \to \alpha$.}\]
The point is to show that the limit $\lim_{r \to \alpha} h(r)$ along rationals $r$ exists.\\
 
We choose $k \in \mathbb{N}$ ``large'' and keep it fixed. We will let $r \to \alpha$ in the end, which by the irrationality of $\alpha$ implies that the length $L$ of the continued fraction expansion of $r$ tends to infinity; thus we can assume for the rest of the argument that with our fixed choice of $k$ we have $k < L$ for all the $r$ we study. We set
\[D = D_k = \max\left\{D(k,\tfrac{1}{k+1},\alpha),D(k,\tfrac{1}{k+1},\alpha'),k\right\}, \] with $D(k,\varepsilon,\alpha)$ being defined as in the statement of Lemma \ref{lem_transf_princ}.
We now consider $\sup_{r \in I_{D}}h(r) - \inf_{r \in I_{D}}h(r)$, where $I_D$ is the interval around $\alpha$ that consists of all numbers $x$ such that the first $D$ partial quotients of $x$ coincide with those of $\alpha$.  By this choice of $D$, we have $k \leq D \leq L$. By Corollary \ref{cor_evil_good_2} and Corollary \ref{cor_r'}, there are suitable numbers $\eta_k$ (which may change from line to line in the following statements, but neither depends on $D$ nor on $L$) such that $\eta_k \to 0$ as $k \to \infty$, and such that 
    \[\sum_{N < q_L}P_N^2(r) \in (1 \pm \eta_k)\sum_{t/2 \leq j < t}\sum_{N \in \mathcal{G}_{j,L}(k,r)}P_N^2(r),\]
and
   \[\sum_{N < q_{L}'}P_N^2(r') \in (1 \pm \eta_k)\sum_{t/2 \leq j < t}\sum_{N \in S(\mathcal{G}_{j,L}(k,r))}P_{N}^2(r').\]
Applying Lemma \ref{lem_decomp} we get 
   \begin{equation*}
\sum_{N < q_L}P_N^2(r) \in 
    (1 \pm \eta_k)\sum_{t/2 \leq j < t}\sum_{N_1 \in \mathcal{G}_{j}^{(1)}(k,r)}P_{N_1}^2(r)\cdot \sum_{N_2 \in \mathcal{G}_{j,L}^{(2)}(k,r)}P_{N_2}^2(r),
   \end{equation*}
   and similarly applying Lemma \ref{split_r'} we get
      \begin{equation*}
      \sum_{N < q_{L}'}P_N^2(r') \in 
    (1 \pm \eta_k)\sum_{t/2 \leq j < t}\sum_{N_1 \in S(\mathcal{G}_{j}^{(1)}(k,r))}P_{N_1}^2(r')\cdot \sum_{N_2 \in S(\mathcal{G}_{j,L}^{(2)}(k,r))}P_{N_2}^2(r').
   \end{equation*}
   We stress once more that  $\mathcal{G}_{j}^{(1)}(k,r)$ does not depend on $L$ since $k < L$.
  By Lemma \ref{lem_transf_princ} and the choice of $D$, we obtain (note that $\mathcal{G}_{j}^{(1)}(k,r) = \mathcal{G}_{j}^{(1)}(k,\alpha)$ since $k < D$ and the first $D$ partial quotients of $r$ and $\alpha$ coincide)
  \begin{equation} \label{mka}
\frac{\sum_{N_1 \in \mathcal{G}_{j}^{(1)}(k,r)}P_{N_1}^2(r)}{\sum_{N_1 \in S(\mathcal{G}_{j}^{(1)}(k,r))}P_{N_1}^2(r')} \in\left(1 \pm O\left(\frac{1}{k}\right)\right)\frac{\sum_{N_1 \in \mathcal{G}_{j}^{(1)}(k,\alpha)}P_{N_1}^2(\alpha)}{\sum_{N_1 \in S(\mathcal{G}_{j}^{(1)}(k,\alpha))}P_{N_1}^2(\alpha')}.
\end{equation}
   By Lemma \ref{lem_tail_conv}, 
   \[\frac{\sum_{t/2 \leq j < t}\sum_{N_2 \in \mathcal{G}_{j,L}^{(2)}(k,r)}P_{N_2}^2(r)}{\sum_{t/2 \leq j < t}\sum_{N_2 \in S(\mathcal{G}_{j,L}^{(2)}(k,r))}P_{N_2}^2(r')} \in 1 \pm \eta_k.\]
Combining the above estimates shows
\[\sup_{r \in I_D}h(r) \leq \eta_k + O\left(\frac{1}{k}\right) + M_k(\alpha)\]
   where
   \[M_k(\alpha) := \log\left(\frac{\sum_{t/2 \leq j < t}\sum_{N_1 \in \mathcal{G}_{j}^{(1)}(k,\alpha)}P_{N_1}^2(\alpha)}{\sum_{t/2 \leq j < t}\sum_{N_1 \in S(\mathcal{G}_{j}^{(1)}(k,\alpha))}P_{N_1}^2(\alpha')}\right).\]
Note how \eqref{mka} was used to make sure that $M_k$ (which captures the impact of the extra partial quotient $a_1$ of $\alpha$, in comparison with $\alpha'$) is a function of $\alpha$, and not a function of $r$.\\
   
   By the same arguments, we obtain
   \[\inf_{r \in I_D}h(r) \geq - \eta_k - O\left(\frac{1}{k}\right) + M_k(\alpha),\]
   which proves 
   \begin{equation} \label{length_con}
\sup_{r \in I_D}h(r)  - \inf_{r \in I_D}h(r) < 2\eta_k + O\left(\frac{1}{k}\right).    
   \end{equation}

   If $(r_n)_{n \in \mathbb{N}}$ is now an arbitrary sequence of rationals with $r_n \to \alpha$, then for each given $k$, there exists $N_0(k)$ with $r_n \in I_{D_k}$ for all $n \geq N_0$. Thus with $n \to \infty$ we can also take $k \to \infty$. Since $(I_D)_{D \in \mathbb{N}}$ is a sequence of nested intervals, also $\left( \left[\inf_{r  \in I_{D_k}} h(r), \sup_{r \in I_{D_k}} h(r)\right] \right)_{k \geq 1}$ is a sequence of nested intervals, whose lengths by \eqref{length_con} converge to zero. Thus there is a unique limiting point of $\left( \left[\inf_{r  \in I_{D_k}} h(r), \sup_{r \in I_{D_k}} h(r)\right] \right)_{k \geq 1}$, and since $r_n \in I_{D_k}$ for all sufficiently large $k$, the limit $\lim_{n \to \infty} h(r_n)$ exists and is finite. Equation \eqref{length_con} also shows that the value of the limit does not depend on the specific sequence $(r_n)_{n \in \mathbb{N}}$, but only on $\alpha$. Thus the limit $\lim_{r \to \alpha} h(r)$ along rationals exists and is finite. 

\subsection{The technical proofs}\label{sec_techn}

For this section, we define
\begin{equation}
\label{alpha_delta_def}\begin{split}
    \delta_k &:= \| q_k \alpha\|,\quad 
\alpha_{k} := [a_{k};a_{k+1},a_{k+2},\ldots], \\
\cev{\alpha}_k &:= [0;a_k,a_{k-1},a_{k-2},\ldots,a_1],\quad \lambda_k:=q_k\delta_k,\end{split} \qquad k \ge 1,
\end{equation}
and we have, for all $k \geq 1$, 
\begin{align}
 \label{alternating_approx}
% q_k\alpha &\equiv (-1)^k\delta_k \pmod 1, \\
\lambda_k &= \frac{1}{\alpha_{k+1}+ \cev{\alpha}_{k}},\\
\label{deltaratio} 
\frac{\delta_{k+2}}{\delta_k}&<\frac{1}{2},\\
%\label{deltarec}
%\delta_{k+1}&=\delta_{k-1}-a_{k+1}\delta_k\le \delta_{k-1}-\delta_k,\\
\label{deltasum} 
\delta_k&=\sum\limits_{t=1}^{\infty}a_{k+2t}\delta_{k+2t-1}.        %\\
%\label{abs_delta}
%\delta_k &= a_{k+2}\delta_{k+1} + \delta_{k+2}.
\end{align}
All these formulas are well-known in Diophantine approximation, for a collection of these (and other related formulas) see e.g. \cite[Section 2]{gay_hau}. 

\subsubsection{Shifted Sudler products}

The proofs of Lemmas \ref{lem_evil_good_1} and \ref{lem_decomp} rely on a decomposition technique which was developed in \cite{grepstad_neum} to solve the Erd\H os--Szekeres problem (answering the problem in the negative by proving that $\liminf_{N \to \infty}P_N(\phi) > 0$ for the Golden Mean $\phi$), and which was brought into a more general and explicit form in recent articles such as \cite{AB_quantum,other_AB,AB3,tech_zaf,grepstadII,hauke_extreme,hauke_density,hauke_badly}. This decomposition of the Sudler product into shifted products related to best approximation denominators was already sketched in the heuristics section around Equation \eqref{factor_heur}. Below, we give a precise statement in the formulation of \cite[Proposition 4]{hauke_badly}.

\begin{prop}\label{prop_shifted}
Let $\alpha$ be a fixed irrational and let $N = \sum_{i=0}^{\ell -1} b_{i}q_i$ be the Ostrowski expansion with respect to $\alpha$ of an integer $N$ in the range $q_{\ell-1} \le N < q_{\ell}$.
For $0 \le i \le n$ and $s \in \mathbb{N}$, using the notation from \eqref{alpha_delta_def} we define
\begin{equation*}
\varepsilon_{i,s}(N) := q_i\left(s\delta_i + \sum_{j=1}^{\ell-i-1} (-1)^jb_{i+j} \delta_{i+j}\right),
\end{equation*}
and 

\[P_{q_i}(\alpha,\varepsilon) := \prod_{n = 1}^{q_i} \left\lvert 2 \sin\Big(\pi\Big(n\alpha + (-1)^i\frac{\varepsilon}{q_i}\Big)\Big)\right\rvert.\]

Then we have

\begin{equation*}
    P_N(\alpha) = \prod_{i=0}^{\ell-1}\prod_{s= 0}^{b_i(N)-1} P_{q_i}\bigl(\alpha,\varepsilon_{i,s}(N)\bigr).
\end{equation*}
\end{prop}

We remark that Propositions \ref{prop_shifted} holds in a perfectly analogous form if one starts with a rational $r$ instead of an irrational $\alpha$. In that case of course one can only consider $N < q_L$, and accordingly $i$ is at most $L-1$. 

\subsubsection{Proof of Lemma \ref{lem_tail_conv}}

We make use of the following statement from \cite[Proposition 4.1]{AB3}:
\begin{prop}\label{tailprop} Let $r = [0;c_1,\ldots,c_L]$, let $1 \le \ell <L$, and let $p_{\ell}'/q_{\ell}'$ denote the convergents of $r'$. Assume the following two conditions:
\begin{enumerate}
\item[(i)] $c_{\ell +1} \le (q_{\ell}')^{1/100}$ or $b_{\ell}(N) \le 0.99 c_{\ell +1}$,
\item[(ii)] $c_{\ell +2} \le (q_{\ell +1}')^{1/100}$ or $b_{\ell +1}(N) \le 0.99 c_{\ell +2}$.
\end{enumerate}
Then for any $N \leq q_L -1$ we have
\[ \prod_{b=0}^{b_{\ell}(N)-1} \frac{P_{q_{\ell}} \left( r, \varepsilon^{(r)}_{b,\ell}(N) \right)}{P_{q_{\ell}'} \left( r',  \varepsilon^{(r')}_{b,\ell-1}(S(N)) \right)} = \exp \left( \mathcal{O} \left( \frac{(c_2+\cdots +c_{\ell})^{3/4}}{(q_{\ell}')^{3/4}} + \frac{\log (c_1+1)}{q_{\ell}'} \right) \right) \]
with a universal implied constant.
\end{prop}

Further, we need \cite[Proposition 3.1]{AB3}. Roughly speaking, the proposition asserts that a Sudler product $P_N(\alpha)$ is particularly large if the Ostrowski digits $b_i(N)$ attain a $5/6$-proportion of their maximal possible value, i.e.\ $b_i(N) \approx \frac{5}{6} c_{i+1}$, for all those $i$ for which the maximal possible value of $b_i$ (namely $c_{i+1}$) is ``large''. For a detailed discussion and the heuristics behind this (and for the connection with the hyperbolic volume of knot complements) we refer to \cite{AB3}. 

\begin{prop}[Local $5/6$-principle]\label{local5/6prop} Let $r = [0;c_1,\ldots,c_L]$. Let $1 \le \ell <L$ be such that $c_{\ell+1} \ge 7$, and set $b_\ell^*:=\lfloor (5/6) c_{\ell+1} \rfloor$. Let $0 \le N <q_L$.
\begin{enumerate}
\item[(i)] If $b_{\ell+1}(N)<c_{\ell+2}$, then $N^*=N+(b_\ell^*-b_\ell(N))q_\ell$ satisfies
\[ \begin{split} \log P_{N^*} &(\alpha ) - \log P_N (\alpha ) \\ & \begin{split} \ge 0.2326 \frac{(b_\ell^*-b_\ell(N))^2}{c_{\ell+1}} -C \bigg( &\frac{|b_\ell^*-b_\ell(N)|}{c_{\ell+1}} \left( 1+\log \max_{1 \le m \le \ell} c_m \right) \\ &+ I_{\{ b_\ell(N) \le 1 \}} I_{\{ b_{\ell+1}(N)>0.99 c_{\ell+2} \}} \log c_{\ell+2} + \frac{1}{q_\ell^2} \bigg) \end{split} \end{split} \]
with a universal constant $C>0$.
\item[(ii)] If $b_{\ell+1}(N)=c_{\ell+2}$, then $N^*=N+b_\ell^* q_\ell - q_{\ell+1}$ satisfies
\[ \begin{split} \log\, &P_{N^*} (\alpha ) - \log P_N (\alpha ) \\ \ge\, &0.1615 c_{\ell+1} -C \left( 1 +\log \max_{1 \le m \le \ell} c_m +\log c_{\ell+2}+I_{\{ c_{\ell+2}=1 \}} I_{\{ b_{\ell+2}(N)>0.99 c_{\ell+3} \}} c_{\ell+3} \right) \end{split} \]
with a universal constant $C>0$.
\end{enumerate}
\end{prop}

\begin{proof}[Proof of Lemma \ref{lem_tail_conv}]
    Note that all elements in $\mathcal{G}_{j,L}^{(2)}(k,r)$ start with at least $(j+1)m \geq k/2$ many zeroes as their first Ostrowski digits. In particular, the first and second Ostrowski digits are both $0$ (recall here that we assumed w.l.o.g.\ that $k \geq 20$). Thus 
    $S_{|\mathcal{G}_{j,L}^{(2)}(k,r)}: \mathcal{G}_{j,L}^{(2)}(k,r)\to S(\mathcal{G}_{j,L}^{(2)}(k,r))$ is well-defined and a bijective map, yielding a one-to-one correspondence between $N_2 \in \mathcal{G}_{j,L}^{(2)}(k,r)$ and $S(N_2) \in S(\mathcal{G}_{j,L}^{(2)}(k,r))$.\\
    
    For fixed $N_2 \in \mathcal{G}_{j,L}^{(2)}(k,r)$, according to Proposition \ref{prop_shifted} the quotient $\frac{P_{N_2}(r)}{P_{S(N_2)}(r')}$ decomposes into factors
    \[\frac{P_{N_2}(r)}{P_{S(N_2)}(r')}
    = \prod_{\ell= (j+1)m}^L\prod_{b=0}^{b_{\ell}(N_2)-1} \frac{P_{q_{\ell}} \left(r,\varepsilon^{(r)}_{b,\ell}(N_2) \right)}{P_{q_{\ell}'} \left( r',  \varepsilon^{(r')}_{b,\ell-1}(S(N_2)) \right)}.
    \]
Following the lines of \cite[Proof of Theorem 5.1]{AB3}, we can use Proposition \ref{local5/6prop} to remove the contribution of those $N_2$ where there exists some $\ell$ in the range $(j+1)m \leq \ell < L-1$ such that the conditions of Proposition \ref{tailprop} are not satisfied. %This procedure, as used in \cite{AB3}, is somewhat analog to our Lemma \ref{lem_evil_good_1} and Corollary \ref{cor_evil_good_2}. The aim now is that the ``evil numbers'' are those where the assumptions of Proposition \ref{tailprop} are not satisfied, i.e. we have unusually large partial quotients $c_{\ell+1}$, and nearly maximal corresponding Ostrowski digits. We now compare such evil $N$ with $N'$ which arises from replacing the corresponding Ostrowski digit of $b_{\ell}(N)$ with $\lfloor (5/6) c_{\ell+1} \rfloor$, with Proposition \ref{local5/6prop} showing that $P_N'(\alpha)$ is larger than $P_N(\alpha)$. This allows to reduce to the case of ``good numbers'' which is the set where the assumptions of Proposition \ref{tailprop} hold true, and we can apply this Proposition to get the desired convergence result.\\
We note that for the application of  Proposition \ref{tailprop} or \ref{local5/6prop} no assumption on the Diophantine properties (such as uniform boundedness of partial quotients) of $r$ is necessary (which indeed we could not guarantee, since we are in the tail of the continued fraction of $r$ which we cannot control by the mere knowledge of $r$ being ``close'' to the badly approximable $\alpha$). Accordingly, the proof from \cite{AB3} carries over verbatim.
Proposition \ref{tailprop} provides exactly what we need, since the arising error terms form a convergent series. After summation over $N_2$ and $j$, this proves the lemma. Note that here we used the fact that the run of consecutive zeros (where we apply the splitting process) appears at an index $\geq k/2$, which thus grows when $k \to \infty$.
\end{proof}

Next, we focus in more detail on the possible perturbations $\varepsilon_{i,s}$ that arise from the decomposition into shifted Sudler products. We define $\varepsilon_i$ to be admissible (with respect to $\alpha$ and $i$, which is omitted if clear from the context) whenever there exist $N \in \mathbb{N}$ and $0 \leq s \leq b_i(N)-1$ such that $\varepsilon_i = \varepsilon_{i,s}(N)$. The following statement provides upper and lower bounds for admissible perturbations:

\begin{prop}\label{max_perturb}
Let $\alpha$ be a fixed irrational. For all $i\in\mathbb{N}$ and every $N \in \mathbb{N}$ we have 
\begin{equation}
\label{eps_min_max}
 -\lambda_i+\lambda_{i,1}\le \varepsilon_{i,s}(N)\le (a_{i+1}-1)\lambda_i+\lambda_{i,1}
\end{equation}
for all $s$ in the range $0 \leq s \leq b_i(N) -1$, where we denote
\begin{equation}
\label{lambda_kj_def}
\lambda_{i,j}=q_i\delta_{i+j}=\frac{q_i}{q_{i+j}}\lambda_{i+j}.
\end{equation}
\end{prop}

\begin{proof}
    This is an immediate consequence of \eqref{deltasum}; for a detailed calculation, see e.g.\ \cite[Proposition 7]{gay_hau}. Furthermore, note that $\alpha$ can also be replaced by a rational $r$ when $i < L$ and $N < q_L$.
\end{proof}

In order to get better control of $P_{q_\ell}(r,\varepsilon)$, we make use of the following approximation. Roughly speaking, the proposition allows us to pass from the shifted Sudler products $P_{q_\ell}(r,\varepsilon)$ to ``limiting'' functions $H_{\ell}(r,\varepsilon)$, which depend on the continued fraction expansion of $r$ in a more direct way.

\begin{prop}\cite[Proposition 7]{hauke_badly}.\label{limit_H}
Let $r = [0;c_1, \dots, c_L]$, and assume for fixed $1 \leq \ell \leq L$ that $\max_{i \leq \ell}c_i \leq M$. Let 
\begin{equation*}
H_\ell(r,\varepsilon) :=
2\pi \lvert \varepsilon + \lambda_\ell \rvert \prod_{n=1}^{\lfloor q_\ell/2\rfloor } h_{n,\ell}(\varepsilon),
\end{equation*}
where
\begin{equation}
\label{hkdef}
h_{n,\ell}(\varepsilon) = h_{n,\ell}(r,\varepsilon) := \Bigg\lvert\Bigg(1 - \lambda_\ell\frac{\left\{n\cev{r}_\ell\right\} - \frac{1}{2}}{n}\Bigg)^2 - \frac{\left(\varepsilon + \frac{\lambda_\ell}{2}\right)^2}{n^2}\Bigg\rvert.
\end{equation}
Further, let $I \subset \mathbb{R}$ be a compact interval. Then we have 
\[P_{q_\ell}(r,\varepsilon) = H_\ell(r,\varepsilon)\left(1 + \mathcal{O}\left(q_\ell^{-2/3}\log^{2/3}q_\ell\right)\right) + \mathcal{O}(q_\ell^{-2}),\quad \forall \varepsilon \in I.\]
The implied constant depends on $M$ and $I$, but neither on $r$ nor on $\ell$.
\end{prop}

\begin{proof}
The above is a slight modification of \cite[Proposition 7]{hauke_badly}, which is only stated there for badly approximable irrationals. However, the only part where this assumption is used is to ensure that $\cev{r}_\ell = [0;c_{\ell},c_{\ell-1},\ldots,c_1]$ has bounded partial quotients, an assumption which is incorporated above by assuming that $\max_{i \leq \ell}c_i \leq M$.
\end{proof}

\begin{prop}\label{hnk_pos}
Let $r = [0;c_1, \dots, c_L]$. Let $\varepsilon_\ell$ be admissible (with respect to $r$ and $\ell \leq L-2$), and assume that $\max_{i \leq \ell+2}c_i \leq M$. Then for all $n \geq 1$, we have
    \[\varepsilon_\ell + \lambda_\ell > \frac{1}{(M + O(1))^2}\]
and
    \begin{equation}\label{eq_hnk_pos}\Bigg(1 - \lambda_\ell \frac{\left\{n\cev{r}_\ell\right\} - \frac{1}{2}}{n}\Bigg)^2 - \frac{\left(\varepsilon_\ell + \frac{\lambda_\ell}{2}\right)^2}{n^2} > \max\left\{1 - O\left(\frac{1}{n}\right), \frac{1}{(M + O(1))^3}\right\},\end{equation}
    where the implied constants are uniform in $r$ and $\ell$.
\end{prop}

\begin{proof}
     Observe that by Proposition \ref{max_perturb}, we have 
$\varepsilon_\ell \geq - \lambda_\ell + \lambda_{\ell,1}$, thus we have (recall \eqref{alternating_approx} and \eqref{lambda_kj_def})
\[\lambda_\ell + \varepsilon_\ell > \lambda_{\ell,1} = \frac{q_{\ell}}{q_{\ell+1}}\frac{1}{r_{\ell+2} + \cev{r}_{\ell+1}}
\geq \frac{1}{M+1}\frac{1}{M + 1},
\]
which proves the first claim.\\
%$\geq - \frac{1}{1 + \frac{1}{M+1}},$
%thus 
%$\lambda_\ell + \varepsilon_\ell > \frac{1}{(M+2)^2}$
%where we used trivial bounds on $q_\ell/q_{\ell+1}$ and $\lambda_\ell$, respectively.
%Thus $\varepsilon + \lambda_\ell$ is bounded away from zero.\\

Since the quantities $\lambda_{\ell},\{n\cev{r}_\ell\},\varepsilon_{\ell}$ on the left-hand side of \eqref{eq_hnk_pos} are all absolutely bounded, it suffices to show \[\Bigg(1 - \lambda_\ell \frac{\left\{n\cev{r}_\ell\right\} - \frac{1}{2}}{n}\Bigg)^2 - \frac{\left(\varepsilon_\ell + \frac{\lambda_\ell}{2}\right)^2}{n^2} > \frac{1}{(M + O(1))^3}.\]

We will show the above only for $n = 1$, since the estimates for larger $n$ are easier and can be treated by trivial estimates. We 
write $x_\ell = 1 - \lambda_\ell\left(\{\cev{r}_\ell\} - 1/2\right), y_\ell = \varepsilon_\ell + \frac{\lambda_\ell}{2}$.
We claim that $x_\ell^2 - y_\ell^2 > \frac{1}{(2M+O(1))^2}$ for all admissible $\varepsilon_\ell$, with the implied constant being absolute. Indeed, using Proposition \ref{max_perturb}, we get 
$y_\ell \leq (c_{\ell+1}-1/2)\lambda_\ell+\lambda_{\ell,1}$, which implies 
\[x_\ell - y_\ell \geq 1 - c_{\ell+1}\lambda_\ell - \lambda_{\ell,1}
\]
by the trivial estimate $\{\cev{r}_\ell\} < 1$. Note that 
\[\lambda_{\ell,1} = q_\ell \delta_{\ell+1} = q_\ell\left(\delta_{\ell-1} - c_{\ell+1}\delta_{\ell}\right) = q_\ell\delta_{\ell-1} - c_{\ell+1}\lambda_\ell,\]
hence
\[\begin{split}x_\ell - y_\ell -1 &\geq - q_\ell\delta_{\ell-1} = -\cev{r}_\ell \lambda_{\ell-1}
= -\cev{r}_\ell \cdot \frac{1}{r_{\ell}+ \cev{r}_{\ell-1}} \geq -1 + \frac{1}{(M+O(1))^2}.
\end{split}
\]
The claim follows now immediately from $x_\ell^2 - y_\ell^2 = (x_\ell - y_\ell)(x_\ell + y_\ell)$ and \[x_\ell + y_\ell > 1 + \varepsilon_{\ell} > 1 - \lambda_{\ell} \geq \frac{1}{M + O(1)}.\]
\end{proof}

\subsubsection{Key lemmas for the proof of Lemma \ref{lem_evil_good_1}} 

We start preparing the proof of Lemma \ref{lem_evil_good_1}. The key technical estimates are the following two lemmas:

\begin{lem}\label{lem_bounded_shifted}
Let $r = [0;c_1, \dots, c_L]$. Let $\ell < L-2$ be given, and assume that $\max_{i \leq \ell+2}c_i \leq M$ and that $\varepsilon_\ell$ is admissible. Then there exists a constant $C = C(M)$ such that
\[\frac{1}{C} \leq H_{\ell}(r,\varepsilon_\ell) \leq C.\]
The constant is uniform in $\ell$ and $r$.
\end{lem}
Further, we need some (essentially Lipschitz) continuity in the argument of perturbation of the Sudler products.

\begin{lem}\label{lem_lip}
Let $r = [0;c_1, \dots, c_L]$. Let $\varepsilon,\varepsilon'$ be admissible with respect to $r$ and $\ell$ and assume that $\max_{i \leq \ell+2} c_i \leq M$. Then there exists a constant $C = C(M)$ such that
    \begin{equation}\label{tail_lip}\left\lvert\frac{P_{q_\ell}(r,\varepsilon)}{P_{q_\ell}(r,\varepsilon')} - 1\right\rvert \leq C \lvert \varepsilon - \varepsilon'\rvert + \mathcal{O}(q_{\ell}^{-2/3}\log^{2/3} q_{\ell}).\end{equation}
Both the constant $C$ and the implied constant are uniform in $r$ and $\ell$.
    
    Furthermore, for fixed $i$, let $I_k(\alpha)$ denote the set of all (rational and irrational) numbers that coincide with $\alpha$ on the first $k$ partial quotients.
    Then for every $\eta > 0$, there exist $\delta,K_0$, such that for 
    $k \geq K_0$,
    \begin{equation}\label{cont_in_perturb}\sup_{r \in I_k(\alpha)}\sup_{\substack{\varepsilon,\varepsilon' \text{admissible w.r.t.\ $i$ and $r$},\\|\varepsilon - \varepsilon'| < \delta}}\left\lvert
    \frac{P_{q_i}(r,\varepsilon)}{P_{q_i}(r,\varepsilon')} - 1\right\rvert < \eta.
    \end{equation}
\end{lem}

We also use the following lemma, which follows from Taylor approximation of the logarithm function.
\begin{lem}[{\cite[Lemma 9]{hauke_extreme}}]
\label{Hauke2021lem}
Let $(x_n)_{A \leq n \leq B}$ be a finite sequence of real numbers that satisfy $|x_n|<\frac{1}{2}$ and $|x_n|<\frac{c}{n}$ for some $c > 0$. Then
\begin{equation*}
\prod\limits_{n=A}^{B}(1-x_n)\ge 1-\Biggl(\Bigl|\sum\limits_{n=A}^B x_n\Bigr|+\frac{c^2}{A-1}\Biggr).
\end{equation*}
\end{lem}

\begin{proof}[Proof of Lemma \ref{lem_bounded_shifted}]
    We will argue similarly to \cite[Proof of Lemma 8]{hauke_badly}. We use Proposition \ref{hnk_pos} in order to remove all absolute values in the definition of $H_\ell$.
  Using $|\frac{1}{2}- \{n \cev{r}_\ell\}|\leq 1/2$, we find the upper bound
\[h_{n,\ell}(r,\varepsilon) \ls
1 + 2q_\ell\delta_\ell \frac{\frac{1}{2}- \{n \cev{r}_\ell\}}{n} + \frac{\varepsilon_{\max}^2 + \varepsilon_{\max} + q_\ell\delta_\ell}{n^2},
\]
where $\varepsilon_{\max}$ is an upper bound for the absolute value of all admissible $\varepsilon_\ell$, which is absolutely bounded. By $\log(1+x) \ls x$, we thus obtain
\begin{equation}\label{H_upperbound}
\begin{split}H_\ell(r,\varepsilon) &\ll_{\varepsilon_{\max}}
 2\pi (1+ q_\ell\delta_\ell)\cdot \exp\left(2q_\ell\delta_\ell\sum_{n = 1}^{\lfloor q_\ell/2\rfloor} \frac{\frac{1}{2}- \{n \cev{r}_\ell\}}{n}\right).
\end{split}
\end{equation}

Employing summation by parts, we get
\[\sum_{n = 1}^{\lfloor q_\ell/2\rfloor} \frac{\frac{1}{2}- \{n \cev{r}_\ell\}}{n}
\leq \frac{S_{\lfloor q_\ell/2\rfloor}(\cev{r}_\ell)}{\lfloor q_\ell/2\rfloor} + 
\sum_{n = 1}^{\lfloor q_\ell/2\rfloor} \frac{S_n(\cev{r}_\ell)}{n^2},
\]
where 
\[S_n(\cev{r}_\ell) := \left\lvert\sum\limits_{u=1}^{n}\frac{1}{2} -\left\{u \cev{r}_\ell) \right\}\right\rvert.\]
From a classical estimate of Ostrowski \cite{ostrowski} (stated there only for irrationals, but also valid for rationals), together with the assumption of the first $\ell+2$ partial quotients of $r$ being bounded by $M$, we have
\begin{equation}
\label{Ostr_est}
S_n(\cev{r}_\ell)\le\frac{3 M \log{n}}{2}.
\end{equation}
Bounding $S_n$ for $n \leq 10$ trivially, we obtain that
\[\sum_{n = 1}^{\lfloor q_\ell/2\rfloor} \frac{\frac{1}{2}- \{n \cev{r}_\ell\}}{n}
\ll_M 1,
\]
which in view of \eqref{H_upperbound} concludes the proof of the upper bound.\\

For the lower bound, we see that by Proposition \ref{hnk_pos}, $h_{n,\ell}(\varepsilon) \geq 1 - \frac{C}{n}, n \geq 1$, for some $C = C(M) > 0$, thus another application of Proposition \ref{hnk_pos} and Lemma \ref{Hauke2021lem} gives for any admissible $\varepsilon$ and any $2 \leq N_0 \leq q_\ell/2$,
\begin{equation}\begin{split}\label{eq_lower_H}H_\ell(r,\varepsilon) &\geq \frac{1}{(M+O(1))^{3N_0}} \prod_{n = N_0+1}^{q_\ell/2}h_{n,\ell}(\varepsilon) 
\\&\geq \frac{1}{(M+O(1))^{3N_0}} \left(1-\Biggl(\Bigl|\sum\limits_{n=N_0+1}^{q_\ell/2} h_{n,\ell}(\varepsilon) - 1\Bigr|+\frac{C}{N_0-1}\Biggr)\right).\end{split}
\end{equation}
Note that
\[h_{n,\ell}(\varepsilon) = 1 - \lambda_\ell\frac{\left\{n\cev{r}_\ell\right\} - \frac{1}{2}}{n} + O(1/n^2)\]
for an absolute implied constant and all admissible (and hence bounded) $\varepsilon$.
Thus by altering $C$ in \eqref{eq_lower_H}, it remains to establish an upper bound for
\[\left\lvert\lambda_\ell\sum_{n + N_0+1}^{q_\ell/2} \frac{\left\{n\cev{r}_\ell\right\} - \frac{1}{2}}{n}\right\rvert.\]
As in the argument for the upper bound in the first half of this proof, by partial summation and employing \eqref{Ostr_est}, we get
\[\left\lvert\sum_{n + N_0+1}^{q_\ell/2} \frac{\left\{n\cev{r}_\ell\right\} - \frac{1}{2}}{n}\right\rvert
\ll_M \sum_{n \geq N_0+1} \frac{\log n}{n^2}.
\]
Since this series converges, choosing $N_0$ sufficiently large shows in combination with \eqref{eq_lower_H} that 
$H_\ell(r,\varepsilon) \gg_{M} 1$. We remark that the implied constant depends only on $M$, and neither on $r$ nor on $\ell$.
\end{proof}

\begin{proof}[Proof of Lemma \ref{lem_lip}]
First, we need to prove that $P_{q_{\ell}}(r,\varepsilon)$ is uniformly bounded from above, as well as uniformly bounded away from $0$: For sufficiently large $\ell$, this follows from Lemma \ref{lem_bounded_shifted} and Proposition \ref{hnk_pos}. 
For the first finitely many $\ell$, we trivially obtain the upper bound $2^{q_{\ell}} \leq 2^{(M+1)^{\ell}}$, which is uniform in $r$. For the lower bound, we use \cite[Lemma 3]{AB_quantum}, which shows for admissible $\varepsilon$, 
that $P_{q_{\ell}}(r,\varepsilon) \gg_{q_\ell,M} 1$. Note that while the results and proofs in \cite{AB_quantum} are only formulated for quadratic irrationals, the arguments applied there actually only use of the fact that $\max_{1 \leq i \leq \ell+1} a_{i}$ is bounded. Further, since $q_{\ell} \leq (M+1)^{\ell}$, the bound is uniform for all $r$ with $\max_{1 \leq i \leq \ell+1} c_{i} \leq M$.
Using this, we see that it is sufficient to show that
%Since by Lemma \ref{lem_bounded_shifted}, $P_{q_\ell}(r,%\varepsilon)$ is bounded from above and below, it suffices to prove %that there exists $C = C(M)$ such that

\[\lvert \log P_{q_\ell}(r,\varepsilon) - \log P_{q_\ell}(r,\varepsilon')\rvert \leq C \lvert \varepsilon - \varepsilon'\rvert
+ \mathcal{O}(q_{\ell}^{-2/3}\log^{2/3} q_{\ell})
.\]
This follows from $|e^x - e^y| < |x-y|e^{z}$ for some $z \in [x,y]$ by the Intermediate Value Theorem, which is applied to $x = 0, y = \log\left(\frac{P_{q_\ell}(r,\varepsilon)}{P_{q_\ell}(r,\varepsilon')}\right)$.\\

Applying Proposition \ref{limit_H}, we can replace $P_{q_\ell}$ by $H_\ell$ since we may assume $\ell$ to be sufficiently large (since otherwise, according to the discussion above, the statement holds trivially when choosing a sufficiently large constant). Using Proposition \ref{lem_bounded_shifted}, we can exchange the additive error term $\mathcal{O}(q_{\ell}^{-2})$ from Proposition \ref{limit_H} with a multiplicative term $1 + \mathcal{O}(q_{\ell}^{-2})$, which is now absorbed by the second error term. Thus, it remains to prove that

\[\lvert\log H_{\ell}(r,\varepsilon) - \log H_{\ell}(r,\varepsilon')\rvert \leq C |\varepsilon - \varepsilon'|.\]

After removing absolute values by applying Proposition \ref{hnk_pos}, we get
\[\log H_{\ell}(r,\varepsilon) - \log H_{\ell}(r,\varepsilon')
= \log\left(\frac{\varepsilon+\frac{\lambda_\ell}{2}}{\varepsilon'+\frac{\lambda_\ell}{2}}\right)
+ \sum_{n = 1}^{\lfloor q_\ell/2\rfloor}\log\left(\frac{\Bigg(1 - \lambda_\ell\frac{\left\{n\cev{r}_\ell\right\} - \frac{1}{2}}{n}\Bigg)^2 - \frac{\left(\varepsilon + \frac{\lambda_\ell}{2}\right)^2}{n^2}}{\Bigg(1 - \lambda_\ell\frac{\left\{n\cev{r}_\ell\right\} - \frac{1}{2}}{n}\Bigg)^2 - \frac{\left(\varepsilon' + \frac{\lambda_\ell}{2}\right)^2}{n^2}}\right).
\]
By Proposition \ref{hnk_pos}, all numerators and denominators in the formula above are bounded from above as well as bounded away from $0$, with the actual size of the bound depending only on $M$. Thus 

\[\frac{\varepsilon+\frac{\lambda_\ell}{2}}{\varepsilon'+\frac{\lambda_\ell}{2}}
= 1 + \frac{\varepsilon' - \varepsilon}{\varepsilon'+\frac{\lambda_\ell}{2}} = 1 + O_M(\lvert \varepsilon-\varepsilon'\rvert),
\]
and hence
\[\log\left(\frac{\varepsilon+\frac{\lambda_\ell}{2}}{\varepsilon'+\frac{\lambda_\ell}{2}}\right) \ll_M \lvert \varepsilon-\varepsilon'\rvert.\]
Similarly,

\[\frac{\Bigg(1 - \lambda_\ell\frac{\left\{n\cev{r}_\ell\right\} - \frac{1}{2}}{n}\Bigg)^2 - \frac{\left(\varepsilon + \frac{\lambda_\ell}{2}\right)^2}{n^2}}{\Bigg(1 - \lambda_\ell\frac{\left\{n\cev{r}_\ell\right\} - \frac{1}{2}}{n}\Bigg)^2 - \frac{\left(\varepsilon' + \frac{\lambda_\ell}{2}\right)^2}{n^2}}
= 1 + \frac{\frac{\left(\varepsilon' + \frac{\lambda_\ell}{2}\right)^2 - \left(\varepsilon + \frac{\lambda_\ell}{2}\right)^2}{n^2}}{\Bigg(1 - \lambda_\ell\frac{\left\{n\cev{r}_\ell\right\} - \frac{1}{2}}{n}\Bigg)^2}
= 1 + O_M\left(\frac{\lvert \varepsilon-\varepsilon'\rvert}{n^2}\right).
\]
Thus 
\[\log H_{\ell}(r,\varepsilon) - \log H_{\ell}(r,\varepsilon') \ll_M
\sum_{n \geq 1} \frac{\lvert \varepsilon-\varepsilon'\rvert}{n^2}
\ll \lvert \varepsilon-\varepsilon'\rvert.
\]
By exchanging the roles of $\varepsilon,\varepsilon'$, \eqref{tail_lip} follows.\\

To show \eqref{cont_in_perturb}, recall that for admissible $\varepsilon$, we have
$P_{q_{\ell}}(\alpha,\varepsilon) \geq C = C(q_{\ell},\alpha)$. We note that in a neighborhood around irrationals, the admissible range of perturbations changes in a continuous way. In other words, for every $\delta > 0$, we have for $r$ sufficiently close to $\alpha$ that if $\varepsilon$ is admissible with respect to $r$ and $\ell$, then there exists $\varepsilon'$ that is admissible with respect to $\alpha$ with the property that 
$\lvert \varepsilon'-\varepsilon\rvert < \delta$; this follows immediately from the definitions of $\varepsilon_{i,s}(N)$, and the fact that $r$ and $\alpha$ coincide on sufficiently many partial quotients. 
Using continuity arguments, we thus obtain for sufficiently small $\delta,$
\[P_{q_{\ell}}(r,\varepsilon) =
\underbrace{P_{q_{\ell}}(r,\varepsilon) - P_{q_{\ell}}(\alpha,\varepsilon)}_{|.| < \delta}
+ \underbrace{P_{q_{\ell}}(\alpha,\varepsilon) - P_{q_{\ell}}(\alpha,\varepsilon')}_{|.| < \delta} + P_{q_{\ell}}(\alpha,\varepsilon) 
\geq C - 2\delta > C/2 > 0.\]

Hence we are \textit{uniformly} bounded away from $0$ in a neighborhood around $\alpha$. By using the same continuity arguments again, \eqref{cont_in_perturb} follows.

\end{proof}

\subsubsection{Proof of Lemma \ref{lem_evil_good_1}}

By Proposition \ref{prop_shifted}, we have
\begin{equation*}
    P_N(r) = \prod_{i=0}^{L}\prod_{s= 0}^{b_i(N)-1} P_{q_i}\bigl(r,\varepsilon_{i,s}(N)\bigr),
\end{equation*}
as well as
\begin{equation*}
    P_{\pi_j(N)}(r) = \prod_{i=0}^{L}\prod_{s= 0}^{b_i(\pi_j(N))-1} P_{q_i}\bigl(r,\varepsilon_{i,s}(\pi_j(N))\bigr).
\end{equation*}
Note that from the definition of $\varepsilon_{i,s}$, i.e. 
\begin{equation*}
\varepsilon_{i,s}(N)  = q_i\left(s\delta_i + \sum_{j=1}^{n-i} (-1)^jb_{i+j}(N) \delta_{i+j}\right),
\end{equation*}
we see that $\varepsilon_{i,s}(N)$ only depends on the digits $b_j, j > i$. Since the digits of $\pi_j(N)$ only may differ at indices $mj,\ldots,m(j+1)-1$ from the ones of $N$, we have
\[\frac{P_N(r)}{P_{\pi_j(N)}(r)}
= \prod_{i=0}^{m(j+1)-1}\frac{\prod_{s= 0}^{b_i(N)-1} P_{q_i}\bigl(r,\varepsilon_{i,s}(N)\bigr)}{\prod_{s= 0}^{b_i(\pi_j(N))-1} P_{q_i}\bigl(r,\varepsilon_{i,s}(\pi_j(N))\bigr)}.
\]
By the definition of $\pi_j(N)$, this implies
\[\begin{split}&\prod_{i=0}^{m(j+1)-1}\frac{\prod_{s= 0}^{b_i(N)-1} P_{q_i}\bigl(r,\varepsilon_{i,s}(N)\bigr)}{\prod_{s= 0}^{b_i(\pi_j(N))-1} P_{q_i}\bigl(r,\varepsilon_{i,s}(\pi_j(N))\bigr)}
\\= &\prod_{i=0}^{mj-1}\frac{\prod_{s= 0}^{b_i(N)-1} P_{q_i}\bigl(r,\varepsilon_{i,s}(N)\bigr)}{\prod_{s= 0}^{b_i(N)-1} P_{q_i}\bigl(r,\varepsilon_{i,s}(\pi_j(N))\bigr)}\cdot 
\prod_{i=mj}^{m(j+1)-1}\prod_{s= 0}^{b_i(N)-1} P_{q_i}\bigl(r,\varepsilon_{i,s}(N)\bigr).
\end{split}\]
Note that $m(j+1)-1 \leq k$ so we have $c_i \leq M$ for all $i \leq m(j+1)+1$, thus we can apply all auxiliary statements from above.\\

An application of Proposition \ref{limit_H} and Lemma \ref{lem_bounded_shifted} thus proves (note that $mj \geq k/2$ implies that this can be chosen sufficiently large)
\[\prod_{i=mj}^{m(j+1)-1}\prod_{s= 0}^{b_i(N)-1} P_{q_i}\bigl(r,\varepsilon_{i,s}(N)\bigr)
\ll \prod_{i=mj}^{m(j+1)-1} C(M)^{c_{i+1}}
\ll \left(C(M)^{M}\right)^{m} = (\tilde{C}(M))^{m}.
\]

An application of Lemma \ref{lem_lip} proves for $0 \leq i \leq mj-1$ that

\[\frac{\prod_{s= 0}^{b_i(N)-1} P_{q_i}\left(r,\varepsilon_{i,s}(N)\right)}{\prod_{s= 0}^{b_i(\pi_j(N))-1} P_{q_i}\bigl(r,\varepsilon_{i,s}(\pi_j(N))\bigr)}
\leq \prod_{s= 0}^{b_i(N)-1} \left(1 + C(M)\lvert \varepsilon_{i,s}(N) - \varepsilon_{i,s}(\pi_j(N))\rvert+ \mathcal{O}(q_i^{-1/2}) \right).
\]

Using \eqref{deltaratio} and \eqref{deltasum}, we get
\[\lvert\varepsilon_{i,s}(N) - \varepsilon_{i,s}(\pi_j(N))\rvert
\ll q_i\delta_{mj} \ll \left(\frac{1}{\sqrt{2}}\right)^{mj - i},\] since all Ostrowski digits $b_{\ell}$ of
$N, \pi_j(N)$ with 
$\ell < mj$ coincide, and further 
$q_i^{-1/2} \leq \frac{1}{\sqrt{2}^{\lfloor i/2\rfloor}}$.
Since this provides a convergent series, we have 
\[\prod_{i=0}^{mj-1}\frac{\prod_{s= 0}^{b_i(N)-1} P_{q_i}\bigl(r,\varepsilon_{i,s}(N)\bigr)}{\prod_{s= 0}^{b_i(N)-1} P_{q_i}\bigl(r,\varepsilon_{i,s}(\pi_j(N))\bigr)}
\ll_M 1,
\]
which completes the proof of Lemma \ref{lem_evil_good_1}.

\begin{proof}[Proof of Corollary \ref{cor_r'}]
    Note that $S(\pi_j(N)) = \pi_j'(S(N))$ where 
    $\pi_j' := \psi^{-1}_{r'} \circ p_j' \circ \psi_{r'}$ and where
\[\begin{split} p_j': \mathcal{A}_{L-1}(r') &\to \mathcal{A}_{L-1}(r')\\
(b_0,b_1,\ldots,b_{L-2}) &\mapsto (b_0,\ldots,b_{jm -2},0,\ldots,0,b_{m(j+1)-1},\ldots,b_{L-2}).\end{split}\]
Analogously to Lemma \ref{lem_evil_good_1}, we can establish for $N \in \mathcal{E}_L(r)$
\[P_{S(N)}(r') \ll ({O_M(1)})^{m} P_{\pi_j'(S(N))}(r)
= ({O_M(1)})^{m} P_{S(\pi_j(N))}(r')
,\]
Thus we get, as in the proof of Corollary \ref{cor_evil_good_2},
    \[\begin{split} t \sum_{N \in S(\mathcal{E}_L(r))}  P_{N}^2(r') &\ll {({O_M(1)})}^m 
    \sum_{t/2 \leq j <t} \sum_{N \in S(\pi_j(\mathcal{E}_L(r)))}P_{N}^2(r')
    \\&\ll {({O_M(1)})}^{m} \sum_{t/2 \leq j < t} \sum_{N \in S(\mathcal{G}_{j,L}(r))}P_{S(N)}^2(r') \\&\ll {({O_M(1)})}^{m}\sum_{N < q_L'}P_N^2(r'),\end{split}\]
    and we can conclude as in Corollary \ref{cor_evil_good_2}.
\end{proof}

\subsubsection{Proof of Lemmas \ref{lem_decomp} and \ref{split_r'}}

We fix $0 < \eta < 1$ arbitrary, and show that there exists a $K_0$ such that for $k \geq K_0$, 
we have for all $j$ in the range $t/2 \leq j < t$ that
\begin{equation}\label{lem_decomp_reform}   \frac{\sum_{N \in \mathcal{G}_{j,L}}P_N^2(r)}{\sum_{N_1 \in \mathcal{G}_{j}^{(1)}}P_{N_1}^2(r) \cdot \sum_{N_2 \in \mathcal{G}_{j,L}^{(2)}}P_{N_2}^2(r)} \in (1 - \eta, 1 + \eta).\end{equation}
This is clearly equivalent to the statement of Lemma \ref{lem_decomp}.
Writing $N = \sum_{0 \leq i \leq L}b_iq_i$ in its Ostrowski expansion with respect to $r$, we 
define \[N_1 = \sum_{0 \leq i \leq jm}b_iq_i,\quad N_2 = \sum_{(j+1)m \leq i \leq L}b_iq_i,\] and observe that for $N \in \mathcal{G}_{j,L}$, we have
$N = N_1 + N_2$. A variant of 
Proposition \ref{prop_shifted} (described below) yields 
\begin{equation} \label{double_1}
P_N(r) = P_{N_2}(r)\cdot\prod_{i=0}^{jm}\prod_{s= 0}^{b_i-1} P_{q_i}\bigl(r,\varepsilon_{i,s}(N)\bigr);
\end{equation}
to obtain this, we do not factorize $P_N(r)$ over all $i$, but only over those indices $i$ that are of size at most $jm$, so that the product $P_{N_2}(r)$ remains intact. This is possible since the shifts $\varepsilon_{i,s}$ for $i \geq j(m+1)$ depend only on the digits with index exceeding $j(m+1)$, but not on the digits with smaller index (which are captured by $N_1$, but do not play a role for $N_2$). Furthermore, applying Proposition \ref{prop_shifted} to $N_1$ we get
\begin{equation} \label{double_2}
P_{N_1}(r) = \prod_{i=0}^{jm}\prod_{s= 0}^{b_i-1} P_{q_i}\bigl(r,\varepsilon_{i,s}(N_1)\bigr).
\end{equation}
The point here is that the double product in \eqref{double_1} is not the same as the one in \eqref{double_2}, since the Ostrowski digits of $N_2$ contribute to $\varepsilon_{i,s}(N)$ but not to $\varepsilon_{i,s}(N_1)$, and so it is not exactly true that $P_N(r) = P_{N_1}(r) P_{N_2}(r)$; however, since the Ostrowski digits of $N_2$ are separated from those of $N_1$ by a long run of zeros (and zeros do not contribute to the value of $\varepsilon_{i,s}$), it will turn out that $\varepsilon_{i,s}(N) \approx \varepsilon_{i,s}(N_1)$, and that accordingly $P_N(r) \approx P_{N_1}(r) P_{N_2}(r)$.\\

More precisely, by \eqref{deltaratio} and \eqref{deltasum}, we obtain that for all $0 \leq i \leq jm$ and for all $s$, we have
$\lvert\varepsilon_{i,s}(N_1) - \varepsilon_{i,s}(N)\rvert \leq (\sqrt{2})^{-m} (\sqrt{2})^{i - jm}$, as a consequence of the fact that the first Ostrowski digit where $N$ and $N_1$ do not coincide has index at least $m(j+1)$. 
In particular, $\lvert\varepsilon_{i,s}(N_1) - \varepsilon_{i,s}(N)\rvert \to 0$, thus an application of Lemma \ref{lem_lip} in the variant of \eqref{cont_in_perturb} implies that for every fixed $i_0$,
there exists $K_0$ such that for $k \geq K_0$,

\[
\frac{\prod_{i=0}^{i_0}\prod_{s= 0}^{b_i-1} P_{q_i}\bigl(r,\varepsilon_{i,s}(N)\bigr)}{\prod_{i=0}^{i_0}\prod_{s= 0}^{b_i-1} P_{q_i}\bigl(r,\varepsilon_{i,s}(N_1)\bigr)}
 \in 1 \pm \frac{\eta}{2},
\]
which holds uniformly for $r$ that coincide with $\alpha$ on the first $k$ partial quotients.
Thus we are left to treat only ``sufficiently large'' values of $i$. Here we apply 
Lemma \ref{lem_lip} in the variant of \eqref{tail_lip} to deduce
\[\begin{split}\frac{\prod_{i=i_0+1}^{jm}\prod_{s= 0}^{b_i-1} P_{q_i}\bigl(r,\varepsilon_{i,s}(N)\bigr)}{\prod_{i=i_0+1}^{jm}\prod_{s= 0}^{b_i-1} P_{q_i}\bigl(r,\varepsilon_{i,s}(N_1)\bigr)}
&\leq \prod_{i=i_0+1}^{jm} \left(1 + C(M)(\sqrt{2})^{-m} (\sqrt{2})^{i - jm} + \mathcal{O}(q_i^{-1/2})\right)^M
\\&\leq \exp \left((\sqrt{2})^{-m}\mathcal{O}_M(1)\sum_{i=0}^{jm}(\sqrt{2})^{i - jm} + \mathcal{O}_M(q_{N_0}^{-1/2})\right)
\\&= \exp \left(O_M((\sqrt{2})^{-m}) + \mathcal{O}_M(q_{i_0}^{-1/2})\right).\end{split}
\]

Since $M$ is fixed and $m = \lfloor \log \log k\rfloor$, the first term becomes arbitrarily small as $k \to \infty$. Further, the second term can be chosen arbitrarily small by increasing the value of $i_0$ accordingly.
The same holds with the roles of $N_1$ and $N$ exchanged, proving that for $N \in \mathcal{G}_{j,L}^{(1)}$,
\[  \frac{P_N(r)}{P_{N_1}(r) \cdot P_{N_2}(r)} \in (1 - \eta, 1 + \eta).\]

Now, recall from \eqref{bij_theta} that 
\[\begin{split}\theta_j: \mathcal{G}_{j,L}^{(1)} \times \mathcal{G}_{j,L}^{(2)} &\to \mathcal{G}_{j,L}\\
(N_1,N_2) &\mapsto N_1 + N_2
\end{split}\]
is a bijection. This completes the proof of Lemma \ref{lem_decomp} (after renaming $\eta$). The proof of Lemma \ref{split_r'} works precisely along the same lines.

\section*{Acknowledgments}

CA was supported by the Austrian Science Fund (FWF), projects 10.55776/I4945, 10.55776/I5554, 10.55776/P34763, 10.55776/P35322, and 10.55776/PAT5120424. MH was supported by the Austrian Science Fund (FWF), project 10.55776/ESP5134624.

\bibliography{Zagier}
\bibliographystyle{abbrv}

\end{document}